\documentclass[a4paper]{amsart}
\usepackage{amsfonts,amssymb}
\usepackage{url}
\usepackage{accents}
\usepackage{txfonts} % \precneqq
\usepackage[dvips]{graphicx}
\usepackage{color}

\newtheorem{Thm}{Theorem}[section]
\newtheorem{Lem}[Thm]{Lemma}
\newtheorem{Cor}[Thm]{Corollary}
\newtheorem{Facts}[Thm]{Facts}
\newtheorem{Fact}[Thm]{Fact}
\theoremstyle{definition}
\newtheorem{Def}[Thm]{Definition}

\newtheorem*{uAsm}{Assumption}
\newtheorem{Exm}[Thm]{Example}
\newtheorem{Rem}[Thm]{Remark}
\newtheorem*{uRem}{Remark}
\newtheorem*{uNote}{Note}
\newtheorem*{uNotes}{Notes}
\newtheorem{Notes}[Thm]{Notes}

\newcommand{\on}{\mathord{\restriction}}
\newcommand{\bR}{\mathbb{R}}

\newcommand{\bN}{\mathbb{N}}
\newcommand{\bQ}{\mathbb{Q}}

\newcommand{\fnto}{\rightarrow}
\DeclareMathOperator{\maxlength}{maxlength}               % defbereich domain
\DeclareMathOperator{\cov}{cov}               % defbereich domain
\DeclareMathOperator{\add}{add}               % defbereich domain
\DeclareMathOperator{\minlength}{minlength}               % defbereich domain
\DeclareMathOperator{\dom}{dom}               % defbereich domain
\newcommand{\card}[1]{\rvert#1\lvert} 

\newcommand{\al}[1]{{\aleph_{#1}}} 
\newcommand{\om}[1]{{\omega_{#1}}} 

\newcommand{\esm}{\prec}

\newcommand{\ho}{^{\omega}}
\newcommand{\hko}{^{<\omega}}                % hoch <omega

\newcommand{\mS}[2]{\underaccent{\tilde}{\Sigma}^#1_#2}
\newcommand{\mP}[2]{\underaccent{\tilde}{\Pi}^#1_#2}

\newcommand{\DEFEQ}{\coloneqq}
\newcommand{\EQDEF}{\eqqcolon}

\newcommand{\forc}{\Vdash}
\newcommand{\incomp}{\perp}
\newcommand{\comp}{\parallel}
\newcommand{\n}[1]{\underaccent{\tilde}{#1}}
\newcommand{\cplsf}{<\hspace{-0.5ex}\cdot}

\newcommand{\Tmax}{T_\text{max}}  
\newcommand{\Tcldn}{T_\text{cldn}}

\newcommand{\af}{\mathfrak f}  % approximations
\newcommand{\ag}{\mathfrak g}
\newcommand{\ah}{\mathfrak h} 
\newcommand{\myI}{\mathbb{I}^{\mathfrak c}}
\newcommand{\strongmyI}{\mathbb{I}}

\DeclareMathOperator{\cf}{cf}

\DeclareMathOperator{\nwflim}{nwf-lim}
\DeclareMathOperator{\length}{length} 
\DeclareMathOperator{\SUCC}{succ}    

\newcommand{\pB}[2]{B^{#1}_{#2}} % used in the definition of the nwf limit
\newcommand{\pP}{\text{Pos}}
\newcommand{\pD}[2]{\text{Dom}^{#1}_{#2}}

\begin{document}

\date{November 2004}

\title{Saccharinity}
\author{Jakob Kellner}
\address{Kurt G\"odel Research Center for Mathematical Logic\\
Universit\"at Wien\\
W\"ahringer Strasse 25\\
1090 Wien, Austria}
\email{kellner@fsmat.at}
\urladdr{http://www.logic.univie.ac.at/$\sim$kellner}

\author{Saharon Shelah}
\address{Einstein Institute of Mathematics\\
Edmond J. Safra Campus, Givat Ram\\
The Hebrew University of Jerusalem\\
Jerusalem, 91904, Israel\\
and
Department of Mathematics\\
Rutgers University\\
New Brunswick, NJ 08854, USA}%{4}
\urladdr{http://shelah.logic.at/}
\thanks{This paper is dedicated to the memory of Greg Hjorth.
Both authors gratefully acknowledge partial support by the National
Science Foundation Grant No. 0600940.  The first author is partially supported
by by European Union FP7 grant PERG02-GA-2207-224747 and the FWF Austrian
Science Fund project P21651-N13.  The second author is supported by the United
States-Israel Binational Science Foundation (Grant no. 2006108), publication
859.}

\begin{abstract}
We present a method to iterate finitely splitting lim-sup tree forcings along
non-wellfounded linear orders. As an application, we introduce a new method to
force (weak) measurability of all definable sets with respect to a certain
(non-ccc) ideal.
\end{abstract}

\maketitle

\section*{Introduction}

\subsection*{Non-wellfounded iterations}
We introduce a method to iterate lim-sup finitely splitting tree forcings
along linear, non-wellfounded orders.

There is quite some literature about non-wellfounded iteration. E.g., Jech and
Groszek~\cite{MR946221} investigated the wellfounded but non-linear iteration
of Sacks forcings. Building on this, Kanovei~\cite{MR1777770} and
Groszek~\cite{MR1295981} develop non-wellfounded iterations of Sacks forcing. In
spirit, their construction is close to the construction of this paper, 
but the  implementation is quite different.
Zapletal gives an illfounded iteration construction in the framework of
``idealized forcing''~\cite{MR2391923}, it seems that his results give some of
the properties of our construction (e.g., $\omega^\omega$-bounding) for a more
general class of forcings, cf.\ his Theorem~5.4.12.\footnote{note that 
Zapletal's basic construction can be applied to countable orders only,
for longer orders additional work is required, see Section~5.5 there.}
Regarding finite support, Brendle~\cite{MR1928381} developed finite-support
non-wellfounded iteration constructions, based on the second author's method of
iterations along smooth templates~\cite{MR2096454}. Brendle's paper also
contains the important observation by Hjorth (answering a question of Hechler)
that it is impossible to have an illfounded iteration of forcings that all add
dominating reals.

\subsection*{Measurability}
As an application of our method, we introduce a new way to 
force measurability of definable sets.

In the seminal paper~\cite{MR0265151} Solovay proved that in the Levy model 
(after collapsing an inaccessible) every definable set is measurable and has
the Baire property.

In~\cite{MR768264} the second author showed that the inaccessible is necessary
for measurability, but the Baire property of every definable set can be
obtained by a forcing $P$ without the use of an inaccessible (i.e.,  in ZFC).
This forcing $P$ is constructed by amalgamation of universally meager forcings
$Q$.  So  every $Q$ adds a co-meager set of generics and has many
automorphisms, and the forcing $P$ has similar properties to the Levy collapse.
The property of $Q$ that implies that $Q$ can be amalgamated is called
``sweetness'' (a strong version of ccc).
One can ask about other ccc ideals than Lebesgue-null and meager (or their
defining forcings, random and Cohen), and classify such ideals (respectively
forcings) according to whether they behave like Cohen or like randoms see,
e.g., Sweet {\&} Sour~\cite{MR2076408}.

For (non-ccc) ideals corresponding to tree forcings $Q$, forcing measurability
can be much simpler, see Section~\ref{sec:cohen} about the Cohen model. In this
model, all definable set are $Q$-measurable (e.g., Marczewski measurable for
$Q$ = Sacks forcing). The proof is a simpler version of Solovay's: Cohen
forcing is homogeneous and adds subtrees $S\in Q\cap V[G]$ to all $T\in Q\cap
V$ such that all branches of $S$ are Cohen reals.

In this paper, we introduce a new construction that gives a variant of
measurability (weak measurability, as defined in~\ref{def:measurable}) for all
definable sets: Instead of iterating basic forcings $Q$ that have many
automorphisms and add a measure 1 set of generics, we use a $Q$ that adds only
a null set of generics (a single one in our case, and this real remains the
only generic over $V$ even in the final limit). So $Q$ has to be very
non-homogeneous. Instead of having many automorphisms in $Q$, we assume that
the skeleton of the iteration has many automorphisms (so in particular a
non-wellfounded iteration has to be used).  

We use the word Saccharinity for this concept: a construction that achieves the
same effect as (an amalgamation of) sweet forcings, but using entirely
different means.

\subsection*{Acknowledgments}
We thank the referee for pointing out many typos and unclarities, and for
providing section~\ref{sec:cohen}.

\subsection*{Annotated contents}\nopagebreak
\begin{list}{}{\setlength{\leftmargin}{0.5cm}\addtolength{\leftmargin}{\labelwidth}}
\item[Section~\ref{sec:Q}, p. \pageref{sec:Q}:] We define a class of 
finitely splitting tree
forcings with ``lim-sup norm'':  The forcing conditions are subtrees of a basic finitely splitting
tree that satisfy ``along every branch, many nodes have many successors''. 
\item[Section~\ref{sec:nwi}, p. \pageref{sec:nwi}:] We introduce a general
construction to iterate such lim-sup tree-forcings along non-wellfounded total
orders.  It turns out that the limit is proper, $\omega\ho$-bounding and has
other nice properties similar to the properties of the lim-sup tree-forcings
itself.
\item[Section~\ref{sec:myI}, p. \pageref{sec:myI}:]
We define (with respect to a lim-sup tree-forcing $Q$) the ideals $\strongmyI$
and $\myI$ (the $<2^\al0$-closure of $\strongmyI$). These ideals will generally
not be ccc.  We define what we mean by ``$X$ is weakly measurable'' and
formulate our application:
Assuming CH and a Ramsey property for $Q$ (see Section~\ref{sec:tree}), we can
force that all definable sets are weakly measurable.  (This section requires
only Section~\ref{sec:Q}.)
\item[Section~\ref{sec:order}, p. \pageref{sec:order}:]
Assuming CH, we construct an order
$I$ which has
many automorphisms and a 
cofinal sequence $(j_\alpha)_{\alpha\in\om2}$.
We show that the non-wellfounded iteration
of $Q$ along the order $I$ forces that $2^\al0=\al2$, that
$\myI$ is nontrivial, that for every definable set
$X$ ``locally'' either
all or none of the generic reals $\n\eta_{j_\delta}$ are in $X$
and that the set $\{\n\eta_{j_\delta}:\, \delta\in\om2\}$ is of weak measure
1 in the set $\{\n \eta_i:\, i\in I\}$.
\item[Section~\ref{sec:tree}, p. \pageref{sec:tree}:]
We assume a certain Ramsey property for $Q$.
We show that $\{\n \eta_i:\, i\in I\}$
is of weak measure 1. Together with
the result of the previous section this proves the application.
\item[Section~\ref{sec:cohen}, p. \pageref{sec:cohen}:]
We give a brief comparison with the Cohen model.
(This section requires only Sections~\ref{sec:Q}
and~\ref{sec:myI}.)
\end{list}

\section{finitely splitting lim-sup tree-forcings}\label{sec:Q}

We will define a class of finitely splitting tree forcings with ``lim-sup
norm''. The simplest example is Sacks forcing.
Such forcings (and
generalizations) have been investigated by many authors,
e.g.\ in~\cite{MR1613600} under the name
$\bQ^\text{tree}_0$ (see Definition 1.3.5 there).

\subsection{Basics}
\begin{figure}[tb]
\begin{center}
\scalebox{0.4}{\input{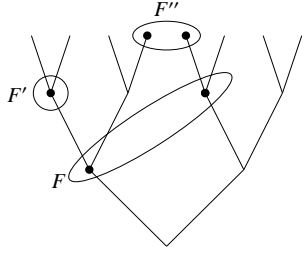}}
\end{center}
\caption{\label{fig:fronts} $F'$ is stronger than $F$, $F''$ is purely stronger
than $F$.}
\end{figure}

Let us first introduce some notation:

\begin{Def}
Let $T\subseteq \omega\hko$ be a tree
(i.e.,  $T$ is closed under initial segments), let
$s,t\in \omega\hko$, $A \subseteq T$.
\begin{itemize}
\item We write sequences as $\langle a_1,\dots,a_n \rangle$ or as
  $(a_1,\dots,a_n)$.
  In particular,  $\langle \rangle$ denotes the empty sequence.
\item $s\preceq t$ means that $s$ is a restriction of $t$ (or equivalently that $s\subseteq t$).
\item $t$ is immediate successor of $s$
	if $t\succeq s$ and $\length(t)=\length(s)+1$.
\item $\SUCC_T(t)$ is the set of immediate successors of $t$ in a tree $T$.
      If the tree $T$ is clear from the context  we will also
      write  $\SUCC(t)$.
\item Nodes $s$ and $t$ are compatible ($s\comp t$), if 
      they are comparable, i.e., if $s\preceq t$ or $t\preceq s$.
        Otherwise, $s$ and
	$t$ are incompatible ($s\incomp t$).
\item The order in forcing notions is usually chosen such that $q<p$
      means that $q$ is stronger than $p$.
      We try to stick to Goldstern's alphabetic
      convention~\cite[1.2]{MR1601976}: Whenever two conditions are comparable
      the notation is chosen so that the variable used for the stronger
      condition comes ``lexicographically'' later.
\item Two forcing conditions $p$ and $p'$ are compatible ($p\comp p'$),
      if there is a $q$ stronger than both $p$ and $p'$. Otherwise,
      $p$ and $p'$ are incompatible ($p\incomp p'$).
\item $T^{[t]}\DEFEQ \{s\in T:\, s\comp t\}$.\quad (So $T^{[t]}$ is a tree.)
      If $T$ is clear, we might also just write $[t]$.
\item $T\on n\DEFEQ\{t\in T:\, \length(t)<n\}$.
\item $A\subseteq T$ is a chain if $s\comp t$ for all $s, t\in A$. 
\item $b\subseteq T$ is a branch if it is a maximal
	chain.\\If there exists a $t\in b$ with length $n$ then
        this $t$ is unique and denoted by $b\on n$.
\item $A\subseteq T$ is an antichain if $s\incomp t$ for all
	$s\neq t\in A$. Unless noted otherwise,
	we will assume that antichains are nonempty.
\item $A\subseteq T$ is a front if it is an antichain
	and every branch $b$ meets $A$ (i.e.,  $\card{b\cap A}=1$).
\item $t\preceq A$ stands for: ``$t\preceq s$ for some $s\in A$''.
\item $\Tcldn^A\DEFEQ\{t\in \omega\hko:\, t\preceq A\}$.\\
        (We will use this downwards-closure 
        only for finite sets $A$. Then $\Tcldn^A$ is a finite tree.)
\item If $A$ and $A'$ are antichains, then $A'$ is stronger than 
	$A$ if for each $t\in A'$ there is a $s\in A$ such that
	$s\preceq t$ (cf.\ Figure~\ref{fig:fronts}).
\item If $A$ and $A'$ are antichains then $A'$ is purely stronger than 
	$A$ if it is stronger and 
	for each $s\in A$ there is a $t\in A'$ such that
	$s\preceq t$ (cf.\ Figure~\ref{fig:fronts}).
\item $\lim(T)$ are the maximal branches of $T$. We use this notation only
      for $T$ that are ``pruned'', i.e., have no finite maximal
      branches; then $\lim(T)\subseteq \omega^\omega$ is the closed set
      corresponding to $T$.
\end{itemize}
\end{Def}

We are only interested in finitely splitting trees 
(i.e.,   $\SUCC(t)$ is finite for all $t\in T$).
Then all fronts are finite.
Note that being a front is stronger than being a maximal antichain.
For example, $\{0^n1:\, n\in\omega\}$ is a maximal antichain in
$2\hko$, but not a front.

\begin{uAsm}
Assume $\Tmax$ and  $\mu$ satisfy the following:
\begin{itemize}
\item $\Tmax$ is a finitely splitting tree.
\item $\mu$ assigns a non-negative real to every subset
	of $\SUCC_{\Tmax}(t)$ for every $t\in \Tmax$.
\item $\mu$ is monotone, i.e.,  if $A\subseteq B$ then $\mu(A)\leq \mu(B)$.
\item The measure of singletons is smaller than 1, i.e.,  $\mu(\{s\})<1$.
\item For all branches $b$ in $\Tmax$, $\limsup_{n\rightarrow \infty}
(\mu(\SUCC(b\on n)))=\infty$.
\end{itemize}
\end{uAsm}

Note that such a $\Tmax$ has to be perfect. 
\begin{Def}\label{def:Q} (The tree forcing $Q$.)
\begin{itemize}
\item
If $T$ is a subtree of $\Tmax$ and $t\in T$,
then $\mu_T(t)$ is defined as the measure of the $T$-successors of $t$,
i.e.,  $\mu_T(t)\DEFEQ\mu(\SUCC_T(t))$.
\item
$Q$ consists of all subtrees $T$ of $\Tmax$ (ordered by inclusion) 
such that along every branch  $b$ of $T$
  \[
    \limsup(\mu_T(b\on n ))=\infty.
  \]
\end{itemize}
\end{Def}
So $\Tmax$ itself is the weakest element of $Q$.

For example, Sacks forcing can be defined in this way:
Set $\Tmax\DEFEQ 2\hko$, and for $t\in \Tmax$ and $A\subseteq \SUCC(t)$ set 
\[
  \mu(A)\DEFEQ\begin{cases}\length(t)&\text{if
}\card{A}=2,\\0&\text{otherwise.}\end{cases}
\]
Then a subtree $T$ of $2\hko$ is in $Q$ iff $T$ is a Sacks tree,
i.e., iff along every branch there
are infinitely many splitting nodes.\footnote{This example is ``atomic''
in the following sense: For a node $s\in T$ there is an $A\subset \SUCC(s)$ such
that $\mu(A)$ is large but $\mu(B)<1$ for every $B\subsetneq A$.
In this paper, we will be interested in ``finer'' norms. In particular
we will require the Ramsey property defined in \ref{def:ramseyprop}.}

\begin{Def}
A (finite or infinite) subtree $T$ of $\Tmax$ is $n$-dense if 
there is a front $F$ in $T$ such that 
$\mu_T(t)>n$ for every $t\in F$.
\end{Def}

\begin{Lem}
\begin{enumerate}
\item A subtree $T$ of $\Tmax$ is in $Q$ iff $T$ is $n$-dense 
for every $n\in\bN$.
\item ``$T\in Q$'' and ``$T\leq_Q S$'' are Borel statements, and ``$S\incomp T $'' is $\mP11$
(in the real parameters $\Tmax$ and $\mu$).
\end{enumerate}
\end{Lem}
\begin{proof}
(1) $\rightarrow$: If $D_n\DEFEQ\{s\in T:\, \mu_T(s)>n\}$ meets 
every branch, then 
\[F_n\DEFEQ\{s\in D_n:\, (\forall s'\precneqq s)\, s'\notin D_n\}\]
is a front.\\
$\leftarrow$:
If $b$ is a branch, then $b$ meets every $F_n$, i.e.,  $\mu_T(b\on m)>n$
 for some $m$.
Since $\mu_T(b\on m)$ is finite,
$\limsup(\mu_T(b \on n))$ has to be infinite.

(2)
Since $\Tmax$ is finitely splitting, ``$F$ is a front'' is equivalent to
``$F$ is a finite maximal antichain''.
\end{proof}

A  finite antichain $A$ can be seen as an approximation to a tree:
``$A$ approximates $T$'' means that $A$ is a front in $T$.
If $A'$ is purely stronger than $A$, then $A'$ gives more 
information about the tree $T$ that is approximated
(i.e.,  every tree approximated by $A'$ is also approximated by $A$). 
And, informally, a stronger antichain approximates smaller
(i.e., stronger) trees.

We will usually identify a finite antichain $F$ and the corresponding
finite tree $\Tcldn^F$.
\begin{Def}\label{def:nsplitting}
\begin{itemize}
\item
A finite antichain $F$ is $n$-dense
if $\Tcldn^F$ is $n$-dense.
\item $\bar F=(F_n)_{n\in \omega}$ is a front-sequence, if
$F_{n+1}$ is $n$-dense and purely stronger than $F_n$.
\item A front-sequence $\bar F$ and a tree $T\in Q$ correspond to each other
if $F_n$ is a front in $T$ for all $n$.
\end{itemize}
\end{Def}

\begin{Facts}
\begin{itemize}
\item
If $F$ is $n$-dense and $F'$ is purely stronger than $F$,
then $F'$ is $n$-dense as well. (This is not true if
$F'$ is just stronger than $F$.)
\item 
If $T\in Q$ then there is a front-sequence corresponding to $T$.
\item
 If $\bar F$ is a front-sequence
then there exists exactly one $T\in Q$ corresponding to $\bar F$,
which we call $\lim(\bar F)$. It is the tree
\[\lim(\bar F)\DEFEQ\{t\in \Tmax:\, (\exists i\in\omega)\, t\preceq F_i\},\]
or equivalently
\[\lim(\bar F)\DEFEQ\{t\in \Tmax:\, (\forall i\in\omega)\, (\exists s\in F_i)\, t\comp s\}.\]
\end{itemize}
\end{Facts}
 
\begin{Lem}\label{lem:Qisproper}
Assume that $Q$ is a
finitely splitting lim-sup tree-forcing. 
\begin{enumerate}
\item
If $T\in Q$ and $t\in T$ then $T^{[t]}\in Q$.
(Sometimes this fact is formulated as ``$Q$ is strongly arboreal''.)
\item\label{item:basicQunion}
The finite union of elements of $Q$ is in $Q$.\footnote{$Q$ is generally
not closed under countable unions.}
\item The generic filter on $Q$ is determined by a real
 $\n \eta$ defined by  $\forc_Q\{\n \eta\}=\bigcap_{T\in G_Q}\lim(T)$;
 or equivalently: $\n\eta$ is the union of the stems of the trees in $G_Q$.
 \\
 It is forced that $\n \eta\notin V$ and that $T\in G_Q$ iff 
 $\n\eta\in\lim(T)$.
 \\ 
 For every $T\in Q$ and $t\in T$, $t\prec \n\eta$
 is compatible with $T$. (In other words:
 $T\not\forc t\not\prec\n\eta$.)
\item (Fusion)
If $(T_i)_{i\in \omega}$ is a decreasing sequence in $Q$ and
$\bar F$ is a front-sequence such that 
$F_i$ is a front in $T_i$ for all $i$,
then $\lim(\bar F)\leq_Q T_i$.
\item (Pure decision) 
If $D\subseteq Q$ is dense, $T\in Q$ and $F$ is a front of $T$,
then there is an $S\leq T$ such that $F$ is a front of $S$ 
and for every $t\in F$, $S^{[t]}\in D$. 
\item $Q$ is proper\footnote{there even are generic conditions for arbitrary
countable transitive ZFC models $M$, similarly to Suslin proper. Sometimes
this is called ``totally proper''.}
and $\omega\ho$-bounding.
\end{enumerate}
\end{Lem}

\begin{proof}[Sketch of proof]
(1) and (2) and (4) are clear. (1) and (2) imply (5).

(3): Let $G$ be $Q$-generic over $V$, and define 
$X\DEFEQ \bigcap_{T\in G}\lim(T)$.
Since $\lim(\Tmax)$ is compact, it satisfies the finite 
intersection property.  So $X$ is nonempty.
For every $T\in G$ and $n\in\omega$ there is exactly one 
$t\in T$ of length $n$ such that $T^{[t]}\in G$.
So $X$ has at most one element.

If $r\in V$, then the set of trees $S\in Q^V$ such that
$r\notin \lim(S)$ is dense: If $r$ is a branch of 
$T\in Q$ then pick an $m$ such that $\mu_T(r \on m)>2$.
Since singletons have measure less than 1,
$r\on m$ has at least two immediate successors in $T$,
and one of them (we call it $t$) is not an initial segment
of $r$. So $S\DEFEQ T^{[t]}$ forces that $\n\eta\neq r$.

Assume towards a contradiction that $\n\eta\in\lim(T)$ for
some $T\in Q^V\setminus G$. Then this is forced by some $S\in G$.
In particular $S$ can not be a subtree of $T$.
So pick an $s\in S\setminus T$. Then 
$S^{[s]}\leq S$ forces that $\n\eta\notin\lim(T)$, a contradiction.

If $T\in Q$ and $t\in T$ then $T^{[t]}$ forces that $t\prec \eta$.

(4) and (5) imply that
$Q$ is $\omega\ho$-bounding and satisfies a version of Axiom A 
(with fronts as indices instead of natural numbers).%
\footnote{In the formulation of fusion and pure decision
we could use the classical Axiom
A version as well: Define
$F^T_n$ to be  the minimal $n$-dense front, i.e., 
\[F^T_n\DEFEQ\{t\in T:\, \mu_T(t)>n\,\&\, (\forall s\precneqq t)\,\mu_T(s)\leq n\},\]
and define  $T\leq_n S$ by $T\leq S$ and $F^T_n=F^S_n$.  It
should be clear how to formulate fusion and pure decision for this notion, and
that this proves Axiom A for $Q$. But in \ref{lem:Qisproper},
we do not use this notion, instead we (implicitly) 
use the following one:
$T\leq_A S$ means that $T\leq S$ and $A$ that 
is a front in both $T$ and $S$.  The reason is that this is the notion that will be
generalized for the non-wellfounded iteration.}
So we get properness.
(We will prove a more general case
in \ref{thm:fusioncor}.)
\end{proof}

So a front can be seen as a finite set of (pairwise incompatible) possibilities
for initial segments of the generic real $\n \eta$.  In the next section we
will generalize this to finite sequences of generic reals instead of a single one. 

\subsection{Some additional facts needed later}
\begin{Lem}\label{lem:treewithoutoldbranch}
If $S\in Q$ and the forcing $R$ adds a new real $\n r\in 2\ho$,
then %there is an $R$-name $\n T$ such that $R$ forces the following:\\
$R$ forces that
there is a $\n T\leq_Q S$ such that
$\lim(\n T)\cap V=\emptyset$, and moreover $\lim(\n T)\cap V$
remains empty in every extension of the universe.
\end{Lem}
\begin{proof}
Assume $S$ corresponds to the front-sequence $\bar F$.
Without loss of generality we can assume that 
along every branch in $S$ there is exactly one split 
between $F_{n-1}$ and $F_n$
and this split has measure $>n$.

\begin{figure}[tb]
\begin{center}
\scalebox{0.7}{\input{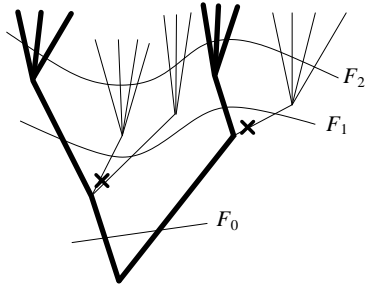}}
\end{center}
\caption{\label{fig:tree1} An example for $S$ and its subtree $\protect\n T$ (bold) when $\protect\n r(0)=0$.}
\end{figure}

We define an $R$-name of a sequence of finite antichains $(\n F'_n)$ 
the following way (cf.\ Figure~\ref{fig:tree1}):
If $n$ is even, we ``take all splits'', i.e., 
 $\n F'_n$ is the set of nodes in $F_n$ that 
are compatible with $\n F'_{n-1}$.
If $n$ is odd, then we add no splittings at all:
for every $s\in \n F'_{n-1}$ we put exactly one successor $t\in F_n$  of $s$ 
into $\n F'_n$, namely the one continuing the $\n r(\frac{n-1}{2})$-th
successor of the (unique) splitting node over $s$.
It is clear that the sequence $(\n F'_n)$ defines a subtree $\n T$ of $S$
that is in $Q$.

Assume $V'$ is an arbitrary extension of $V$ containing an
$R$-generic filter $G$ over $V$. If
$\eta\in \lim(\n T[G])\cap V$, then
$\n r[G]$ can be decoded in $V$ using $S$ and $\eta$. This is a contradiction
to $\forc_R \n r\notin V$.
\end{proof}

We will also need the following family of definable dense subsets
of $Q$:

\begin{Def}\label{def:spl}
Fix a recursive bijection $\psi$ from $\omega$ to $2\hko$. 
Assume that $f:\omega\fnto\omega$ is strictly increasing and 
that $A\subseteq
\omega$. 
\begin{itemize}
\item For $g\in 2\ho$, define $A^{\psi}_g\DEFEQ \{n\in\omega:\, \psi(n)\prec g\}$.
\item  $Q^f_A$ is the set of all $T\in Q$ such that  for all splitting nodes $t\in T$,
  $\length(t)$ is in the interval $[f(n),f(n+1)-1]$ for some $n\in A$.
\item $T\in Q$ has full splitting with respect to $f$
  if for all $n\in\omega$ and $t\in T$ of length $f(n+1)$ there
  is an $s\preceq t$ of length at least $f(n)$ such that
  $\mu_T(s)>n$.
\item 
	$D^\text{spl}_f$ is the set of all $T\in Q$ such that
	either $T\in Q^f_{A^{\psi}_g}$ for some $g\in 2\ho$
	or $T \incomp_Q S$ for all $g\in 2\ho$ and $S\in Q^f_{A^{\psi}_g}$.
\end{itemize}
\end{Def}

Of course the notions $Q^f_{A}$ and $D^\text{spl}_f$ depend
on the forcing $Q$ (i.e.,  on $\Tmax$ and $\mu$), so maybe it
would be more exact to write $Q^f_{A}[\Tmax,\mu]$ etc. However,
we always assume that the $Q$ is understood. 
\ref{def:measurable}).

\begin{Lem}\label{lem:dspl}
Assume that $f:\omega\fnto\omega$ is strictly increasing and  $A,B\subseteq \omega$.
\begin{enumerate}
\item If $g\neq g'$, then $A^{\psi}_g\cap A^{\psi}_{g'}$ is finite.
\item\label{item:dspl1}
        $Q^f_\omega=Q$.
        \quad 
	If $A$ is finite then $Q^f_A=\emptyset$.
\item\label{item:dspl2}
	$Q^f_A\cap Q^f_B= Q^f_{A \cap B}$.\quad
	If $A\subseteq B$, then $Q^f_A\subseteq Q^f_B$.
\item\label{item:dspl3}
	If $T\leq_Q S$ and $S\in Q^f_A$ then $T\in Q^f_A$.
\item\label{item:dspl4}
	For every $T\in Q$ there is a strictly increasing
        $f$ such that $T$ has full splitting with respect to $f$.
\item\label{item:dspl5}
	If $T\in Q$ has full splitting with respect to $f$
	and $\card{A}=\al0$ then there is an $S\leq_Q T$
	such that $S\in Q^f_A$.
\item\label{item:dspl6}
	$D^\text{spl}_f$ is an (absolute definition of an) open dense subset of $Q$
	(using the parameters $f$, $\Tmax$ and $\mu$).\footnote{$X_0\DEFEQ 
\{A^{\psi}_g:\, g\in 2\ho\}$
is an almost disjoint family, but not
maximal.  So of course $Q(X_0)\DEFEQ \bigcup_{A\in X_0} Q^f_A\subset Q$ is not
dense.  We add the incompatible conditions to get the dense set
$D^\text{spl}_f$.  One could ask whether $Q(X)$ is dense for a m.a.d. family
$X$. The following holds:
\begin{enumerate}
\item For every $f$ there is a m.a.d. family $X$ such that $Q(X)$
is not dense.
\item (CH) For every $f$ there is a m.a.d. family $X$ 
such that $Q(X)$ is dense.
\end{enumerate}
Proof: Fix $f$. 
A node $s\in\Tmax$ has level $m$ if $f(m)\leq \length(s)<f(m+1)$.  $S\in Q$ has
unique splitting if $S$ has at most one splitting point of level $n$ for all
$n\in\omega$. For every $T\in Q$ there is an $S\leq_Q T$ with unique
splitting.

For~(a), fix a $T\in Q$ with unique splitting.
Set $Y\DEFEQ\{A\in[\omega]^\al0:\, (\forall S\leq_Q T)\,S\notin Q^f_A\}$.
$Y$ is open dense in $([\omega]^\al0,\subseteq)$, therefore there is
a m.a.d. $X\subseteq Y$.

For~(b), list $Q$ as $(T_\alpha)_{\alpha\in\om1}$, and build
$B_\alpha\in[\omega]^\al0$ by induction on $\alpha\in\om1$: Find an $S\leq_Q
T_\alpha$ with unique splitting.  If some $S'\leq_Q S$ is in $Q^f_{B_\beta}$
($\beta<\alpha$) (or equivalently in $Q^f_{\bigcup_{i\in l} B_{\beta_i}}$ 
for some $l\in\omega$, $\beta_0,\dots,\beta_{l-1}<\alpha$),
then just pick any almost disjoint $B_\alpha$.
Otherwise
enumerate $(B_\beta)_{\beta\in\alpha}$ as $(C_n)_{n\in\omega}$,
and construct $B_\alpha$ and $S'\leq_Q S$ inductively:
At stage $n$, add a split of $S$ to 
$S'$ whose level is not in $\bigcup_{m\leq n} C_m$, and use some
bookkeeping to guarantee that $S'\in Q$. Let $B_\alpha$
be the set of splitting-levels of $S'$.}
\item\label{item:dspl7}
	  In any extension $V'$ of $V$ the following holds:
        If $r\in 2\ho\setminus V$ and
        $S\in Q^f_{A^{\psi}_{r}}$, then
	$T\incomp_Q S$ for all $T\in V\cap D^\text{spl}_f$.
\end{enumerate}
\end{Lem}

\begin{proof}
(1)--(4) and (6) are clear.

(5):
Let $T$ be an element of $Q$. 
Assume we already constructed $f(n)$. Let $N$ be the maximum
of $\mu_T(t)$ for $t\in T\on f(n)$.
There is an $N+n+1$-dense front $F$ in $T$.
Let $f(n+1)$ be the maximum of $\{\length(t):\, t\in F\}$.

(7):
``$T$ is incompatible with all $S\in Q^f_{A^{\psi}_g}$'' is absolute, since it is equivalent to
\[(\forall  g\in 2^\omega)\, (\forall S\subseteq \Tmax)\ \left[S\notin Q^f_{A^{\psi}_g}\ \vee\  
   T\incomp_Q S\right],\]
which is a $\mP11$ statement.

(8):
Let $r\in 2\ho\setminus V$
and $T\in V\cap D^\text{spl}_f$. 
If $T\in Q^f_{A^{\psi}_g}$ for some $g\in 2^\omega\cap V$, 
then $g\neq r$, so $A^{\psi}_g\cap A^{\psi}_r$ is finite and
$Q^f_{A^{\psi}_{r}}\cap Q^f_{A^{\psi}_{g}}$ is empty.
If on the other hand $T$ is incompatible with all $S\in Q^f_{A^{\psi}_g}$ in $V$
 then this holds in $V'$ as well.
\end{proof}

Assume $f'(n)\geq f(n)$ for all $n\in\omega$. Define $h(n)$ by induction:
$h(n+1)\DEFEQ f'(h(n)+1)$.
If $T$ has full splitting with respect to $f$, then $T$ has full splitting with respect to $h$:
$h(n)\leq f(h(n))$, since $f$ is strictly increasing. 
$f(h(n)+1)\leq f'(h(n)+1)=h(n+1)$, and there
are $h(n)$-dense splits between 
the levels $f(h(n))$ and $f(h(n)+1)$. So there are $n$-dense splits between the levels 
$h(n)$ and $h(n+1)$.  So we get:

\begin{Lem}\label{lem:fullsplitinV}
If $V'$ is an $\omega\ho$-bounding extension of $V$
and $T\in Q^{V'}$, then there is a strictly increasing $h\in V$ such that
(in $V'$) $T$ has full splitting with respect to $h$.
\end{Lem}

\section{A non-wellfounded Iteration}\label{sec:nwi}

In this section we introduce a general construction to iterate lim-sup
tree-forcings $Q_i$
(as defined in the last section) along non-wellfounded linear orders $I$.  It
turns out that the limit $P$ is proper, $\omega\ho$-bounding and has other nice
properties similar to the properties of $Q_i$ itself. If $I$ is wellfounded,
then $P$ is equivalent to the usual countable support iteration of
(the evaluations of the definitions) $Q_i$.

\subsection{Conditions and approximations, the nw-iteration}
\begin{Def}
  Let $I$ be a linear order.
  For $i\in I$ we set $I_{<i}\DEFEQ \{j\in I:\, j<i\}$ and 
  analogously we define $I_{\leq i}$ and $I_{>i}$.
  We also set $I_{<\infty}\DEFEQ I$.
\end{Def}

For every $i\in I$ we fix a finitely splitting lim-sup tree-forcing $Q_i$
(to be more exact, we fix a pair $\Tmax^i,\mu^i$).
In the application of this paper, each $Q_i$ will be the same forcing $Q$.

\begin{Def}\label{def:precondition}(Pre-condition)
We call $p$ a pre-condition, if $p$ is a function, the domain of $p$
	is a countable\footnote{this includes finite and empty.} subset of $I$,
	and for each $i\in\dom(p)$,
	$p(i)$ consists of the following:
        \begin{itemize}
          \item $\pD p i$, a countable subset of $\dom(p)\cap I_{<i}$, and
	  \item a (definition of a) Borel function $\pB p i: 
            \left(\omega^\omega\right)^{\pD p i}\fnto Q_i$
        \end{itemize}
\end{Def}

\begin{Rem}
The idea is that we calculate the condition $\pB p i\in Q_i$ using countably
many generic reals $(\eta_j)_{j\in\pD p i}$ that have  already been produced at
stage $i$.  The forcing conditions $p$ of the non-wellfounded iterations will
be pre-conditions that satisfy additional properties, in particular: all $\pB p
i$ are continuous (on a certain Borel set), i.e., if we want to know $\pB p i$
up to some finite height we only have to know $(\eta_i\on m)_{i\in u}$
for some finite $u$ and $m\in\omega$.  Moreover, we will assume that we will
have ``wellfounded continuity parts''. This will be explained in the following,
here just an example: Assume that $I=\omega^*=\{\dots,3,2,1,0\}$, and each
$\Tmax^i=2^{<\omega}$.  Let $p$ be the pre-condition with $\pD p n=\{n+1\}$,
i.e., $\pB p n$ only depends on the generic real $\n\eta_{n+1}$, and $\pB p
n(x)=[0]$ if $x(0)=0$ and $\pB p n(x)=[1]$ otherwise.  Then $p$ is continuous,
but will not be a valid condition, since it is not well founded enough.
\end{Rem}

\begin{figure}[tb]
\begin{center}
\scalebox{0.4}{\input{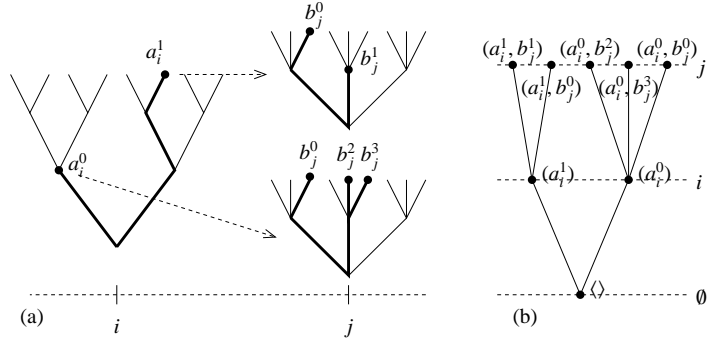}}
\end{center}
\caption{\label{fig:approx_both} An approximation $\ag$:
 $u=\{i,j\}$, $\Tmax^i=2\hko$, $\Tmax^j=3\hko$.
 \protect\\
 $\pP(\ag)=\pP_{\leq j}(\ag)=\{(a_i^1,b_j^0),(a_i^1,b_j^1),(a_i^0,b_j^2),(a_i^0,b_j^3),(a_i^0,b_j^0)\}$.
 \protect\\
 (a): viewed as function:
 $\ag(i)(\langle\rangle)=\{a_i^0,a_i^1\}$, $\ag(j)(\langle a_i^1\rangle)=\{b_j^0,b_j^1\}$ etc.
 \protect\\
 (b): viewed as tree, the heights labeled with $\{\emptyset\}\cup u$.
 }
\end{figure}

We now define finite ``approximations'' to conditions of the iteration; they
will have the same role for the iteration that finite antichains have for $Q$
(see, e.g., Lemma~\ref{lem:Qisproper}).
The following definition looks rather unpleasant, but really is quite simple,
as Figure~\ref{fig:approx_both} hopefully demonstrates.  (We first define
approximations as functions as in~(a) of the figure; sometimes it is more
useful to think of them as trees as in~(b), which will be described
in~\ref{fact:approxtree}.)

\begin{Def}\label{def:approximation} (Approximation)
  \begin{itemize}
    \item
      $\ag$ is an approximation, if 
      $\ag$ is a function with finite domain
      $u\subseteq I$ of the following form:
      Let $i_0$ be the smallest element of $u$. We set
      $\pP_{<i_0}(\ag)\DEFEQ\{\langle\rangle\}$.
      By induction
      on $i\in u$, we assume that $\pP_{<i}(\ag)$ is a 
      set of sequences indexed by the set $\{j\in u: j<i\}$,
      and require the following:
      $\ag(i)$ is a function from $\pP_{<i}(\ag)$
      to finite antichains in $\Tmax^i$, and we set
      \[
        \pP_{\leq i}(\ag)\DEFEQ\{\bar a^\frown b:\,
        \bar a\in \pP_{<i}(\ag),\, b\in \ag(i)(\bar a)\}.
      \]
      If $j$ is the successor of $i$ in $u$, 
      we set $\pP_{<j}(\ag)$ to be $\pP_{\leq i}(\ag)$.
    \item For any $i\in I\cup\{\infty\}$, we define $\pP_{< i}(\ag)$ as
      $\pP_{\leq j}(\ag)$, where $j=\max(\dom(\ag)\cap I_{<i})$
      (or as $\{\langle\rangle\}$, if $\dom(\ag)\cap I_{<i}$ is empty).
      We set $\pP(\ag)\DEFEQ\pP_{<\infty}(\ag)$ and
      call it the set of possibilities of $\ag$.
    \item\label{item:incr_def_of_approx}
      If $i\notin \dom(\ag)$ or $\bar a\notin \pP_{<i}(\ag)$  
      we set $\ag(i)(\bar a)\DEFEQ\{\langle\rangle\}$
      (i.e.,  the front in $\Tmax^i$ consisting only of the root.
       This corresponds to ``no information'').
    \item
      Let $\ag$ be an approximation,
      $J\subset I$, and $\bar \eta=(\eta_i)_{i\in J}$ a
      sequence of reals.  
      Then ``$\bar \eta$ is compatible with 
      $ \ag $'', if there is an $\bar a\in \pP(\ag)$ 
      such that $a_i\prec \eta_i$ for all 
      $i\in \dom(\ag)\cap J$.
      If in addition $J\supseteq \dom(\ag)$, then this $\bar a$ is uniquely defined and called $\bar \eta\on \ag$.
      If $J\supseteq \dom(\ag)\cap I_{<i}$, then $\bar a\on I_{<i}$ is uniquely defined, and therefore we can set 
      $\ag(i)(\bar \eta)\DEFEQ\ag(i)(\bar a\on I_{<i})$.
      
      If $\bar b=(b_i)_{i\in J}$ is a sequence of elements of $\omega\hko$,
      we define $\bar b$ to be compatible with $\ag$  if
      there is a sequence $\bar \eta$ extending $\bar b$ and compatible with $\ag$.
      If $J\supseteq \dom(\ag)$ and additionally 
      each $b_i$ is long enough, then such a $\bar b$ defines a unique
      $\bar a\in \pP(\ag)$ called $\bar b\on \ag$;
      if  $J\supseteq \dom(\ag)\cap I_{<i}$ and additionally 
      each $b_i$ is long enough, then 
      we can define $\ag(i)(\bar b)$ as above.
    \item
      If $\ag$ and $\ag'$ are both approximations, then
      ``$\ag'$ is stronger than $\ag$'' if $\dom(\ag')\supseteq \dom(\ag)$
      and for all 
      $\bar b\in \pP(\ag')$
      there is an $\bar a\in \pP(\ag)$ such that
      $\bar b\succeq \bar a$ (i.e.,  $b_i \succeq a_i$ for all $i\in \dom(\ag)$). 
      In this case $\bar a$ is $\bar b\on \ag$.
      
      Equivalently, $\ag'$ is stronger than $\ag$ iff 
      for all $i\in\dom(\ag)$ and all  
      $\bar b\in\pP_{<i}(\ag')$ there is a (unique) $\bar a\in \pP_{<i}(\ag)$ 
      such that $\bar b\succeq \bar a$ and the antichain
      $\ag'(i)(\bar b)$ is stronger than $\ag(i)(\bar a)$.
    \item
      $\ag'$ is purely stronger than $\ag$ if $\ag'$ is stronger than $\ag$ and
      for all $i\in\dom(\ag)$ 
      and $\bar b\in \pP_{<i }(\ag')$
      the front $\ag'(i)(\bar b)$ is purely stronger
      than $\ag(i)(\bar b\on \ag)$.
    \item 
      For $u\subseteq \dom(\ag)$, 
      $\maxlength_{u}(\ag)$ is 
      $\max(\{\length(a_i):\, i\in u,\bar a\in\pP(\ag)\})$.
      \\
      $\maxlength(\ag)$ is $\maxlength_{\dom(\ag)}(\ag)$.
      Analogously we define $\minlength(\ag)$.
    \item $\ag$ is $n$-dense at 
      $i\in I$, if $i\in \dom(\ag)$
      and for all $\bar a\in \pP_{<i}(\ag)$, $\ag(i)(\bar a)$
      is $n$-dense for $Q_i$. (See Definition~\ref{def:nsplitting}.)
    \item For all
      $\bar a=(a_i)_{i\in u}$ 
      such that $a_i\in \Tmax^i$ there is a (unique) approximation
      $\ag$ such that $\pP(\ag)=\{\bar a\}$. We will call this 
      approximation $\bar a$ as well.
  \end{itemize}
\end{Def}

\begin{Facts}\label{fact:partord}
  \begin{itemize}
    \item ``stronger'' is a partial order on the set of approximations;
      the same holds for ``purely stronger''.
    \item If $\ah$ is stronger than $\ag$, then
      all $\bar \eta$ compatible with $\ah$ are compatible with $\ag$.
  \end{itemize}
\end{Facts}

We could equivalently define approximations as trees, cf.\
Figure~\ref{fig:approx_both}(b): Given an approximation $\ag$, we can define an
approximation-tree with $u=\dom(\ag)$ labeling the heights above the root, 
and the set of
nodes at height $i_n\in u$ is $\pP_{\leq i_n}(\ag)$; the tree order is just
extension of sequences. Every such approximation-tree corresponds to
an approximation:
\begin{Fact}\label{fact:approxtree}
  Consider a finite tree where the heights above the root 
  are labeled by the increasing sequence $u=\{i_1,\dots,i_n\}$
  in $I$. 
  Assume that each node at height $i_m$ is a sequence $(a_j)_{j=i_1,\dots,i_{m}}$ and that the tree order is the extension relation.
  Then this tree corresponds to an approximation, iff each branch has maximal
  height and the successors of each node at level $i_{n-1}$ form an antichain
  in $\Tmax^{i_n}$.
\end{Fact}
In particular, if we take a subset of the (maximal) branches in the
approximation-tree $\ag$, we get a ``sub-approximation'' $\ah$.
A single branch $\bar a$ is a special case of such a sub-approximation.

\begin{Def} (Approximation to $p$)
	Let $p$ be a pre-condition.
  \begin{itemize}
    \item
      $\ag$ approximates $p$, or: $\ag$ is a $p$-approximation,
      if $\dom(\ag)\subseteq \dom(p)$ and $\ag$ is an
      approximation with the following property: If $i\in \dom(\ag)$, $\bar a\in 
      \pP_{<i}(\ag)$, and $\bar\eta=(\eta_j)_{j\in\pD p i}$ is compatible with $\bar a$,
      then $\ag(i)(\bar a)$ is a front in $\pB p i(\bar \eta)$.
    \item 
      $\ag$ is an indirect approximation to $p$ witnessed by $\ag'$,
      if $\ag'$ approximates $p$ and $\ag'$ is purely stronger than $\ag$.
  \end{itemize}
\end{Def}

\begin{Exm}
  The following trivial example should demonstrate the difference  between
  approximation and indirect approximation:
  Assume each $\Tmax^i$ is $2^{<\omega}$, and $p$ is a condition with
  $\dom(p)=\{i,j\}$ for some $i<j$ in $I$.
  Accordingly $\pD p i$ has to be empty, and 
  $\pB p i$ is constant; we set it to have constant value $[1]$.
  We set $\pD p j=\{i\}$ and $\pB p j(x)=[x(0)]$, i.e., if
  the real $x$ starts with $0$ then $\pB p j$ is $[0]$ and otherwise
  it is $[1]$.
  We define the approximation $\ag$ by $\pP(\ag)=\{(\langle\rangle, 1)\}$ and
   $\ah$ by $\pP(\ah)=\{(1, 1)\}$.  
  Then $\ag$ indirectly approximates $p$, witnessed by $\ah$.
\end{Exm}

Now we can define the forcing $P$, the non-wellfounded
countable support limit along $I$:

\begin{Def} (The nwf-iteration $P=\nwflim_I(Q_i)$)
  \begin{itemize}
    \item $p\in P$ means:\\
      $p$ is a pre-condition, and
      for all finite $u\subseteq \dom(p)$, $i\in u$ and $n\in \omega$
      there is a $p$-approximation $\ag$ such that $\dom(\ag)\supseteq u$, $\ag$ is
      $n$-dense for $i$, and $\minlength_{u}(\ag)>n$.
    \item For $p,q\in P$, $q\leq p$ means:\\
      for all $p$-approximations $\ag$
      there is a  $q$-approximation $\ah$ which is stronger than $\ag$
      (so in particular, $\dom (q)\supseteq \dom(p)$).
    \item $q\leq_{\ag} p$ if $q\leq p$ and $\ag$ indirectly approximates $p$
	and $q$.
  \end{itemize}
\end{Def}

\begin{uRem}
  The definition of $q\leq_P p$ is {\em not} equivalent
  to ``for all $i$ and $\bar\eta$,
  $\pB q i(\bar\eta)$ is a subtree of $\pB p i(\bar\eta)$.''
  (Informally speaking, we are only interested in ``the generic $\bar\eta$,
  not in ``all $\bar\eta$''.)
  We will see 
  in Lemma~\ref{lem:muchstuff}(\ref{item:specificstuff2})
  that  $q\leq_P p$ is equivalent to: for each $i\in I$
  it is forced by $q\on P_{<i}$ that $\pB q i(\bar\eta)$
  is a subtree of $\pB p i(\bar\eta)$, where $\bar\eta$
  is the generic sequence up to $i$.
\end{uRem}

\begin{Facts}
  \begin{itemize}
    \item
      $\leq$  is transitive, and for a fixed approximation $\ag$
      the relation $\leq_{\ag}$ is transitive as well. 
    \item If $\ah$ is
      purely stronger than $\ag$ then $\leq_{\ah}$ implies $\leq_{\ag}$.
    \item
      For every $p\in P$, the approximations of $p$ are directed:
      If $\ag$ and $\ag'$ both (indirectly)
      approximate $p$, then there is a $\ah$ approximating $p$ that is
      (purely) stronger than both $\ag$ and $\ag'$. In fact, {\em every}
      $p$-approximation $\ah$ has this property if
      $\dom(\ah)\supseteq \dom(\ag)\cup\dom(\ag')$ and if
      $\minlength_{\dom(\ag)\cup\dom(\ag')}(\ah)$ 
      is large enough.
  \end{itemize} 
\end{Facts} 

So in particular for every $p\in P$ there is an approximating sequence:
\begin{Def}
  An approximating sequence for $p\in P$
  is a sequence $(\ag_n)_{n\in\omega}$ of approximations of $p$
  such that $\ag_{n+1}$ is purely stronger than $\ag_n$, and
  $\ag_{n+1}$ is $n$-dense for each $i\in\dom(\ag_n)$, and 
  $\dom(p)=\bigcup_{n\in\omega} \dom(\ag_n)$.
\end{Def}

An approximating sequence contains all
relevant information about $p$. 
In particular, $\ag$ is an indirect approximation to $p$ iff
there is an $n$ such that $\ag_n$ is purely stronger than $\ag$.
So if $p$ and $q$ both have the approximating sequence $(\ag_n)_{n\in\omega}$,
then $p=^*q$ (i.e., $p\leq q$ and $q\leq p$), furthermore $\ag$
indirectly approximates $p$ iff it indirectly approximates $q$.

Approximating sequences provide an equivalent definition for $P$:
\begin{Def} (Alternative definition of the nwf-iteration $P$)
Define the p.o. $P'$ as follows:
$\bar{\ag}\in P'$ iff $\bar{\ag}$ is a sequence of approximations
$(\ag_n)_{n\in\omega}$ such that
$\ag_{n+1}$ is purely stronger than $\ag_n$ and
$n$-dense for every $i\in\dom(\ag_n)$.
We define $\bar{\ah}<\bar{\ag}$ as: For every $n$ there is an $m$ such that
$\ah_m$ is stronger than $\ag_n$.
\end{Def}

\begin{Lem}\label{lem:sequofapprox}
There is a dense embedding%
\footnote{$\phi$ is even an isomorphism modulo $=^*$, where
$p=^*q$ if $q\leq p$ and $q\leq p$.}
$\phi:\, P'\fnto P$.
\end{Lem}
\begin{proof}
Given a sequence $\bar{\ag}\in P'$, define 
$p=\phi(\bar{\ag})$ the following way:
$\dom(p)=\bigcup \dom(\ag_n)$. For $i\in \dom(p)$,
set $\pD p i\DEFEQ \dom(p)\cap I_{<i}$. 
Define $T= \pB p i(\bar \eta)$ as follows:
If $\bar \eta$ is compatible with all 
$\ag_n$, then let $T$ be
$\{t\in \Tmax^i:\, (\exists n\in\omega)\, t\preceq \ag_n(i)(\bar\eta)\}$. 
Otherwise, let $n$ be maximal such that $\bar \eta$ is compatible
with $\ag_n$,  and let $T$ be $\{t\in \Tmax^i:\,
(\exists s\in \ag_n(i)(\bar \eta))\, t\comp s\}$. Clearly,
$\pB p i$  is
a Borel function, 
$\pB p i(\bar\eta)\in Q_i$
and each $\ag_n$ approximates $p$.
Therefore $(\ag_n)_{n\in\omega}$ is an approximating sequence for $p\in P$.
It is clear that $\phi$ preserves the order.

Let $\psi$ map $p\in P$ to any approximating sequence for $p$.
$\psi: P\fnto P'$ preserves order as well and
$\phi(\psi(p))=^* p$. Therefore $\phi$ is a dense embedding.
\end{proof}

\begin{Notes}\label{notes:conditions}
\begin{enumerate}
\item\label{tmp4_11} If $\ag$ indirectly approximates $p$,
then there is a $q=^*p $ such that $\ag$ approximates $q$.
(Just let $q$ correspond to an approximating sequence of
$p$ starting with $\ag_0=\ag$.) 
\item
It doesn't matter whether the $\ag_n$ in an approximating sequence 
are approximations to $p$ or just indirect approximations.
\item
It doesn't matter whether $\ag_{n+1}$ proves $n$-density
for every $i\in\dom(\ag_n)$ or for just some $i_n$,
provided that the sequence $(i_n)_{n\in\omega}$ covers $\bigcup\dom(\ag_n)$
infinitely often.
\item
In Definition~\ref{def:precondition} 
of pre-condition, instead of requiring $\pB p i$
to be a function into $Q_i$, we could have defined 
$\pB p i$ to be a function to subtrees of
$\Tmax^i$. The additional ``$n$-dense'' requirements on a condition
guarantee $\pB p i(\bar \eta)\in Q_i$
anyway (for generic sequences $\bar \eta$).
\item
Every approximation $\ag$ can be interpreted as 
a condition in $P$, by
\[\pB \ag i(\bar \eta)\DEFEQ \{t:\, t\comp \ag(i)(\bar \eta)\}
	\text{ for }i\in\dom(\ag).\]
(Where we set $\ag(i)(\bar \eta)\DEFEQ \{\langle\rangle\}$
if $\bar \eta$ is incompatible with $\ag$.)
Then $\ag$ approximates itself.
\item\label{tmp3_4}
For any approximation $\ag$ and $u\subseteq I$ finite we can assume $u\subseteq
\dom\ag$: Just set $\ag(i)$ to be the constant function with value
$\{\langle\rangle\}$ for $i\notin\dom\ag$. (Recall that
$\{\langle\rangle\}$ is the ``trivial front'' corresponding
to ``no information''.)
\item\label{ijewqwqtqwt}
If $\ag$ and $\ah$ are approximations, we can assume without
loss of generality that $\dom(\ag)=\dom(\ah)$.
\item\label{ijewqwqtqwt22}
For any $U\subseteq I$ countable and $p\in P$
 we can assume without loss of generality
that $\dom(p)\supseteq U$. This is clear if
$p$ is interpreted as a sequence of Borel-functions:
just set $\pB p i$ to be (the constant function with value)
$\Tmax^i$ for
$i\notin\dom(p)$. 
If $p$ is interpreted as sequence $(\ag_n)_{n\in\omega}$ of approximations, we
have to set $\ag_n(i)$ to be (the constant function with value)
$\Tmax^i \cap \omega^{k(n)}$ for some sufficiently large $k(n)$.  (Using 
$\{\langle \rangle\}$ does not work here, since it does not satisfy $n$-density.)
\item\label{ijewqwqtqwt23}
So if $q\leq p$ we can assume $\dom(q)=\dom(p)$, and if  $p$ is interpreted as
sequence $(\ag_n)_{n\in\omega}$ and $q$ as  $(\ag_n)_{n\in\omega}$ then we can
assume $\dom(\ag)=\dom(\ah)$.
\end{enumerate}
\end{Notes}

\subsection{Fusion and pure decision}

We have seen: Every $p\in P$ corresponds to a 
purely increasing sequence $(\ag_n)$ of approximations such that
$\bigcup\dom(\ag_n)=\dom(p)$ and $\ag_{n+1}$ is $n$-dense for
$\dom(\ag_n)$.
The approximating sequences immediately prove a version of fusion:

\begin{Lem}\label{lem:fusion}
(Fusion) 
Assume that $(p_n)_{n\in\omega}$ 
is a sequence of conditions,
$(\ag_n)_{n\in\omega}$ a sequence of approximations, and $i_n\in\dom(\ag_n)$ 
 such that:
\begin{itemize}
\item $p_{n+1}\leq_{\ag_n} p_n$,
\item $\ag_{n+1}$ is purely stronger than $\ag_n$ and $n$-dense for $i_n$,
\item $(i_n)_{n\in\omega}$ covers $\bigcup \dom(p_n)$ infinitely 
often.
\end{itemize}
Then there is a condition $p_\omega$ such that 
$p_\omega\leq_{\ag_n} p_n$
for all $n$.
\end{Lem}

\begin{proof}
We already know that the sequence $(\ag_n)_{n\in\omega}$ of approximations
defines a condition $p_\omega$ 
such that each $\ag_n$ approximates $p_\omega$.
If $\ah$ approximates $p_n$, then some $\ag_m$
is stronger than $\ah$. Then $\ag_m$ approximates
$p_\omega$, so $p_\omega\leq p_n$. 
\end{proof}

\begin{Def}\label{def:subapprox}
  $\ah$ is sub-approximation of $\ag$ if $\pP(\ah)\subseteq \pP(\ag)$.  (So in
  particular $\dom(\ag)=\dom(\ah)$.)
\end{Def}
Obviously any sub-approximation of $\ag$ is stronger than $\ag$.
In the interpretation of approximations as trees, a sub-approximation is just
a nonempty subset of the (maximal) branches, see Fact~\ref{fact:approxtree}.

\begin{Lem} (Sub-approximation)
  Assume that $\ag$ indirectly approximates $p$ and that $\ah$ is a
  sub-approximation of $\ag$. Then there is 
  a weakest condition stronger than $p$ and approximated by $\ah$,
  which we call $p\on \ah$.
\end{Lem}
\begin{proof}
  Without loss of generality, we can think of $p$ as
  an approximation-sequence $(\ag_n)_{n\in\omega}$ with
  $\ag=\ag_0$. 
  We define approximations $\ah_n$ as follows:
  $\ah_n$ consists of those nodes in the approximation-tree
  $\ag_n$ that are compatible with an element of $\ah$.
  Then $p\on \ah$ is the 
  sequence $(\ah_n)_{n\in\omega}$.
\end{proof}

A special case of a sub-approximation is a singleton:
\begin{Def}\label{def:restrtoa}
  Assume that $\ag$ (indirectly) approximates $p$ and $\bar a\in \pP(\ag)$.  
  We can interpret $\bar a$ as an approximation, a sub-approximation
  of $\ag$. Instead of $p\on \bar a$ we also write
  $p^{[\bar a]}$.
\end{Def}

\begin{Cor}\label{cor:restrtoa}(Specialization and pure decision)
Assume that $\ag$ indirectly approximates $p$ and that $\bar a\in\pP(\ag)$.
\begin{enumerate}
\item $p^{[\bar a]}\in P$, $p^{[\bar a]}\leq p$ and
	$\bar a$ indirectly approximates $p^{[\bar a]}$.
	If $q\leq p$ and $\bar a$ indirectly approximates 
	$q$, then $q\leq p^{[\bar a]}$.
\item If $q\leq_{\ag}p$, then 
	$q^{[\bar a]}\leq p^{[\bar a]}$.
\item\label{item:rbetween} If $q\leq p^{[\bar a]}$
	then there is a $r\leq_{\ag} p$ such that $r^{[\bar a]}=^*q$.
\item\label{item:specialipredense} 
	The set
	$\{p^{[\bar a]}:\, \bar a\in\pP(\ag)\}$ is predense below $p$.
\item\label{item:singletondense} 
	Abusing notation, we denote with $(i,a)$ the approximation 
	$\ag$ with domain $\{i\}$ such that $\ag(i)(\langle\rangle)=\{a\}$.
        For all $i\in I$, $n\in \omega$ the following set is dense:
	\[
          \{p\in P:\, (\exists a\in\omega^n)\,  (i,a)\text{ approximates }p\}.
        \]
	(Or, in the notation introduced later: We can densely determine the generic
        $\n\eta_i$ up to $n$.)
\item\label{item:puredecision} (Pure decision)
	If $D\subseteq P$ is open dense, and $\ag$ indirectly approximates
	$p$, then 
	there is an $r\leq_{\ag} p$ such that 
	$r^{[\bar a]}\in D$ for all $\bar a\in\pP(\ag)$.
\end{enumerate}
\end{Cor}

\begin{proof}
(1) and (2) follow easily from the definition.

(\ref{item:rbetween}) 
We set $r$ to be $q$ ``below $\bar a$'' and $p$ otherwise.  Let $p$ correspond
to $(\ag_n)_{n\in\omega}$ with $\ag_0=\ag$, and $q$ corresponds to
$(\ah_n)_{n\in\omega}$ with $\ah_0=\bar a$ such that each $\ah_n$ is stronger than
$\ag_n$.
According to Note~\ref{notes:conditions}(\ref{ijewqwqtqwt23}), we can assume
that $\dom(\ah_n)=\dom(\ag_n)=u_n$.
We define by induction on $n$ a sub-approximation $\af_n$ of $\ag_n$:
 Let $i_0$ be minimal in $u_n$. So
  $\pP_{<i_0}(\af_n)=\{\langle\rangle\}$. By induction on $i\in u_n$, define for
  all $\bar b\in\pP_{<i}(\af_n)$ 
  \[\af_n(i)(\bar b)\DEFEQ\begin{cases}
  \ag_n(i)(\bar b)&
   \text{if }\bar b\text{ is incompatible with }\ah_n,\\
\ah_n(i)(\bar b) \cup \{t\in \ag_n(i)(\bar b):\, t\incomp
\ah_n(i)(\bar b)\}&
   \text{otherwise.}
\end{cases}\]
It is clear that the possibilities of
$\af_n$ follow $\ah_n$ up to some $i\in\dom \ag_n$
and from then on become incompatible with $\ah_n$ and follow $\ag_n$.
To be more exact: $\bar b\in \pP(\af_n)$ iff $\bar b
\in\pP(\ag_n)$ and
for some $i\in\dom(\ag_n)\cup\{\infty\}$, $\bar a\on I_{<i}$
is in $\pP(\ah_n)$ and either $i=\infty$ or $a_i \incomp
\ah_n(i)(\bar a)$. From this it follows that 
$\af_n$ is purely stronger than $\ag_n$, and that the
$\af_n$ are an approximating sequence
(converging to some $r\leq p$).

(\ref{item:specialipredense}) If $\ag$ indirectly approximates $p$ and $q\leq
p$, then there is a $\ah$ stronger than $\ag$ approximating $q$.  Let $\bar
b\in \pP(\ah)$ and $\bar a=\bar b\on \ag\in\pP(\ag)$. Then $q^{[\bar
b]}\leq q,p^{[\bar a]}$.

(\ref{item:singletondense})
Let $\ah$ approximate $p$ such that $\minlength_{\{i\}}(\ah)>n$.
Let $\bar a\in\pP(\ah)$. Then $(a_i)$ indirectly 
approximates $p^{[\bar a]}\leq p$.
By \ref{notes:conditions}(\ref{tmp4_11}) we can find a
$q=^* p$ such that $(a_i)$ approximates $q$.

(\ref{item:puredecision})
Let $\pP(\ag)=\{{\bar a}_0,\dots,{\bar a}_l \}$.
	Pick $q_0\leq p^{[{\bar a}_0]}$ in $D$, and
	$r_0\leq_{\ag}p$ as in (\ref{item:rbetween}).
	So $r_0^{[{\bar a}_0]}\in D$. 
	Pick $q_1\leq r_0^{[{\bar a}_1]}$ in $D$
	and $r_1\leq_{\ag}r_0$ as above, etc.
	Then $r_l$ has the required property.
\end{proof}

\begin{uRem}
	Similarly, we can define conjunctions of two approximations
	$\ag,\ag'$.  More specifically: let us call $\ag$ and $\ag'$ compatible
	if there is an $\ah$ stronger than both $\ag$ and $\ag'$. Then for
	every compatible pair $\ag,\ag'$ there is a weakest approximation
	$\ag\wedge \ag'$ stronger than $\ag$ and $\ag'$. If $p$ and $q$ have
	incompatible approximations, then they are incompatible (in $Q$).  This
	can be used to define the conjunction of an approximation and a
	condition (if the condition $p$ corresponds to the sequence $\ag_n$,
	let $p\wedge \ah$ correspond to the sequence $\ag_n\wedge \ah$; it is
	the weakest condition stronger than $p$ that is approximated by $\ah$).
	Similarly one can define the conjunction of two conditions. However,
	all of this will not needed in this paper.
\end{uRem}

\subsection{Restrictions}

We now list some trivial properties of $P$ regarding restriction:
\begin{Def} For $i\in I\cup\{\infty\}$ we define $P_{<i}\DEFEQ 
	\{p\in P:\,\dom(p)\subseteq I_{<i}\}$.
	In particular, $P=P_{<\infty}$.
	Analogously we define $P_{\leq i}$ for $i\in I$.
\end{Def}

\begin{Facts}\label{facts:restriction}
(Restriction) Assume $p,q\in P$ and $i,j\in I\cup\{\infty\}$.
\begin{itemize}
\item
If $\dom(q)\supseteq \dom(p)$, $q\on \dom(p)=p$ and 
  $\ag$ approximates $p$, 
	then $q\leq_{\ag} p$.
\item 
$ p\on I_{<i}\in P_{<i}$ and $p\leq  p\on I_{<i}$.
\item If $p'\leq p$ then $p'\on I_{<i}\leq p\on I_{<i}$. 
	If $p\in P_{<i}$ then $p\on I_{<i}=p$.
\item Let $q\in P_{<i}$, $q\leq p\on I_{<i}$. Define 
	$q\wedge p\DEFEQ  q\cup p\on I_{\geq i}$. Then
	$q\wedge p\in P$ is the weakest condition stronger than both
        $q$ and $p$.
\item $p\on I_{<i}$ 
	is a reduction of $p$ (i.e.,  $r'\in P_{<i}$ and $r'\leq p\on
	I_{<i}$ implies $r' \comp p$).
\item In particular,
  $P_{<i}\cplsf P_{<j}$ (i.e.,  $P_{<i}$ is a complete subforcing of $P_{<j}$)
  for $i\leq j$.
\item If $p\on I_{<i} \comp q\on I_{<i}$
	and $\dom(p)\cap \dom(q)\subseteq I_{<i}$, then $p\comp q$.
\item Similar facts hold for $P_{\leq i}$. E.g., if
	$i<j$, then $P_{\leq i}\cplsf P_{<j}$.
\end{itemize}
\end{Facts}

\begin{Def}
Assume that $j\in I\cup\{\infty\}$ and $i<j$, and that $G_{<j}$ is a 
$P_{<j}$-generic filter over  $V$.
\begin{itemize}
\item
Since $P_{<i}$ is a complete 
subforcing of $P_{<j}$, the filter 
$G_{<j}\cap P_{<i}\EQDEF G_{<i}$ is $P_{<i}$-generic over $V$.
We set $V_{<i}\DEFEQ V[G_{<i}]$.
The canonical $Q_i$-generic filter over $V_{<i}$ is called $G(i)$.
Analogously we can define $V_{\leq i}$ and $G_{\leq i}$
(which turns out to be $V_{<i}[G(i)]$ and $G_{<i}\ast G(i)$, respectively).
\item
In $V_{<j}$ or $V_{\leq i}$
we define $\eta_i$ to be the union of all 
$t\in\omega\hko$ such that
$(i,t)$ is an approximation\footnote{as in
Corollary~\ref{cor:restrtoa}(\ref{item:singletondense})}
of $p$ for some $p\in G_{<j}$ (or $G_{\leq i}$). 
\end{itemize}
\end{Def}

\begin{Lem}\label{lem:muchstuff}
Let $i,j,G_{<j}$ be as above, $p\in G_{<j}$, and set $\bar \eta=(\eta_l)_{l<j}$.
\begin{enumerate}
\item $\eta_i$ is a well-defined real. In particular
  we can calculate $\pB q i(\bar \eta\on \pD q i)$ for
  all $q\in P$; abusing notation, we will just write
  $\pB q i(\bar\eta)$.
\item If $\ag$ indirectly approximates $p$, then $\bar \eta$ is compatible with $\ag$.
\item $\{\eta_i\}=\bigcap \{\lim\pB q i(\bar \eta):\, q\in G_{<j},\,i\in\dom(q)\}$.
\item $q\in G_{<j}$ iff 
  $\eta_i\in\lim(\pB q i(\bar \eta))$ for all $i\in\dom(q)$.
\item\label{item:specificstuff} (in $V$): 
  $q\leq_P p$ iff $\dom(q)\supseteq\dom(p)$ and 
  $q\forc \n\eta_i\in\lim(\pB p i(\bar{\n\eta}))$ for all $i\in\dom(p)$.
\item\label{item:specificstuff2} (in $V$): 
  $q\leq_P p$ iff $\dom(q)\supseteq\dom(p)$ and 
  $q\on I_{<i}\forc \pB q i(\bar{\n\eta})\subseteq \pB p i(\bar{\n\eta})$ for all $i\in\dom(p)$.
\end{enumerate}
\end{Lem}

\begin{proof}
(1) By \ref{cor:restrtoa}(\ref{item:singletondense}),
	the set of conditions $q$ such that for some $t$ of length $n$
	the approximation
	$(i,t)$ approximates $q$ is dense.
	Therefore $\eta_i$ is infinite. Also,
	if $s\incomp t$, if $(i,t)$ is an approximation of $q$,
	and if $(i,s)$ is an approximation of $q'$,
	then $q$ and $q'$ are incompatible.
	This shows that $\eta_i$ is indeed a real.

(2) According to \ref{cor:restrtoa}(\ref{item:specialipredense}), the set
	$\{p^{[\bar a]}:\, \bar a\in\pP(\ag)\}$ is 
	predense below $p$. Let $\bar a$ be such that
	$p^{[\bar a]}\in G$. 
	Any $q\in G$ that is stronger than $p^{[\bar a]}$ and
	decides $\n \eta_i$ up to the length of $a_i$ forces  that
	$\n \eta_i\supset a_i$. So $\bar \eta$ is compatible
	with $\bar a$ and therefore with $\ag$.

(3) Let $n\in \omega$. We have to show that
	$\eta_i\on n\in \pB q i(\bar \eta)$.
	First pick an approximation $\ag$ of $q$ 
	with $\minlength_{\{i\}}(\ag)\geq n$.
	We already know that $\bar \eta$ is compatible with $\ag$,
	in particular  $\eta_i$ is compatible with 
	$\ag(i)(\bar \eta)$. 
	And $\ag(i)(\bar \eta)$ is a front in 
	$\pB q i(\bar \eta)$, since $\ag$ approximates $q$.
     It remains to be seen that the intersection on the right-hand side
     is a singleton; this is clear by genericity.

(5) One direction follows immediately from the definition
    of the order in $P$: Assume that $q\leq p$ and that $i\in\dom(p)$.
    Assume towards a contradiction that $r\leq q$
    forces that $\n\eta_i\notin\lim(\pB p i(\bar{\n\eta}))$,
    more specifically that $\n\eta_i\on M
    \notin\pB p i(\bar{\n\eta})$ for some $M$ (already determined
    by $r$). Pick a $p$-approximation $\ag$ 
    that has minimal height 
    greater than $M$ at position $i$;
    and an $r$-approximation $\ah$ stronger than $\ag$.
    Pick $\bar b\in\pP(\ah)$ and let $\bar a\in \pP(\ag)$
    be the restriction. 
    Then $r^{[\bar b]}$ forces that
    $\n\eta_i\on M<a_i<b_i$
    for any $\bar b\in\pP(\ah)$,
    but $a_i\in\ag(i)(\bar a)$ which is a front in $\pB p i(\bar{\n\eta})$,
    a contradiction.

    For the other direction, 
    let $\ag$ approximate $p$ and $\ah$ approximate $q$ such that
    $\dom(\ah)\supseteq \dom(\ag)$ and the length of $\ah$ is sufficiently
    large on $\dom(\ag)$.
    Then $\ah$ must be stronger than $\ag$, which shows that $q\leq p$.

(4) follows from (3) and (5); (6) follows from (5) (see also the proof
of Lemma~\ref{lem:lem35}).
\end{proof}

\subsection{Properness, bounding, continuous reading}
As immediate consequence of fusion and pure decision we get:
\begin{Thm}\label{thm:fusioncor}
\begin{enumerate}
\item\label{item:omegaomegabounding} $P$ is $\omega^\omega$-bounding.
  For every $p$ and $P$-name $\n\tau$ for an $\omega$-sequence of ordinals
  there is a $q\leq p$ such that $q$ reads $\n\tau$ continuously.%
\footnote{
  In more detail: Let $(\n \tau(n))_{n\in \omega}$ be a sequence
  of $P$-names for ordinals and $p\in P$. Then there is a 
  $q\leq p$ corresponding to a sequence $(\ag_n)_{n\in\omega}$
  of approximations, and there are functions $f_n$ from
  $\pP(\ag_n)$ into the ordinals
  such that $q^{[\bar a]}$ forces $\n \tau(n)=f_n(\bar a)$
  for all $\bar a\in \pP(\ag_n)$.
  If each $\n\tau(n)$ is a natural number then this defines
  (in $V$) a continuous function $F$ from $(\omega^\omega)^{\dom(q)}$
  into $\omega^\omega$
  such that $q$ forces that $F(\bar \eta\on \dom(q))=\bar{\n\tau}$.
}

\item\label{item:nonewreals} 
	Assume that the cofinality\footnote{We always mean
	the ``upwards cofinality'', i.e.,  the
	minimal size of an upwards cofinal subset.
	$A\subset I$ is upwards cofinal if for every $i\in I$
	there is an $a\in A$ such that $a\geq i$.}
        of $I$ is $\geq\al1$, that $G$ is $P$-generic over $V$ 
	and that $r\in\bR^{V[G]}$. Then there is an $i\in I$
	such that $r\in\bR^{V_{<i}}$.
\item $P$ is proper.\footnote{$P$ even is non-elementary-proper (nep),
i.e., there are generic conditions for all (non-transitive, non-elementary,
but ord-transitive) countable ZFC models; cf.\ \cite{MR2115943} or\ \cite{nep}.}
\item $P$ forces that $\n\eta_i$ is a $Q_i$-generic real over $V_{<i}$.
\item 
If $I=I_1+I_2$, then $\nwflim_I(Q_i)\cong\nwflim_{I_1}(Q_i)\ast \nwflim_{I_2}(Q_i)$,
the forcing-composition of $\nwflim_{I_1}(Q_i)$ and (the
$\nwflim_{I_1}(Q_i)$-name for) $\nwflim_{I_2}(Q_i)$.
\item\label{item:fusiocorCSI} If 
$I=\Sigma_{\beta\in\epsilon} J_\beta$ is the concatenation of the orders
	$J_\beta$ along the ordinal $\epsilon$,
	then $\nwflim_I(Q_i)$ is equivalent to the countable support limit 
        $(P_\beta,\n Q'_\beta)_{\beta\in\epsilon}$, where
        $\n Q'_\beta$ is (the $P_\beta$-name for)
        $\nwflim_{J_\beta}(Q_i)$.
\item If $I$ is well-founded, then $\nwflim_I(Q_i)$ is the 
	countable support limit of the $Q_i$.
\end{enumerate}
\end{Thm}

\begin{proof}
(1)
Fix for every countable 
subset $J$ of $I$ an enumeration $\{j_m:\,m\in \omega\}$,
and denote $\{j_m:\, m\in n\}$ by $\text{first}(n,J)$.

Assume $\n\tau$ is a name of a real and $p\in P$.
We have to show that there is a $p_\omega\leq p$ and an
$f\in \omega\ho$
such that $p_\omega\forc \n\tau(n)<f(n)$.
Let $p_0\leq p$, $f(0)\in\omega$ 
be such that $p_0\forc \n\tau(0)=f(0)$,
and let $\ag_0$ approximate $p_0$.
Assume that $\ag_n$ and $p_n$ are already defined.
We define $p_{n+1}\leq_{\ag_n} p_n $, $f(n)$ and
$\ag_{n+1}$
the following way: Let $p_{n+1}\leq_{\ag_n} p_n $
be such that
$p^{[\bar a]}_{n+1}$ decides $\n\tau(n)$
for every $\bar a\in \pP(\ag_n)$, 
see \ref{cor:restrtoa}(\ref{item:puredecision}).
Let $f(n)$ be the maximum of the possible values for
$\n\tau(n)$. Let $\ag_{n+1}$ be a
$p_{n+1}$-approximation 
stronger than $\ag_n$ which is $n$-dense at every
$i\in \text{first}(n,\dom(p_1))\cup\dots\cup \text{first}(n,\dom(p_n))$.
Then the sequence $(p_n)_{n\in\omega}$ 
satisfies the conditions for fusion \ref{lem:fusion}
so there is a $p_\omega\leq p_n$. Clearly, $p_\omega\forc \n\tau(n)\leq f(n)$.

The same argument shows continuous reading of $\omega$-sequences: Now we do not
require $\n\tau(n)$ to be a natural number, and we do not care about the
maximum possible value; the rest is the same.

(2) The $p_\omega$ above completely determines $\n\tau$,
	so if $p_\omega\in P_{<i}$, then $p_\omega \forc_P \n\tau\in V_{<i}$.

(3) is very similar to the above: Assume that $N\esm H(\chi)$ and
	$p_0\in N$. 
	Let $\{D_m:\, m\in\omega\}$ enumerate the dense sets in $N$.
	Assume $p_n$, $\ag_n\in N$ are already defined.
	Find (in $N$) $p_{n+1}\leq_{\ag_n} p_n$ such that
 $p^{[\bar a]}_{n+1}\in D_n$ for all $\bar a\in\pP(\ag_n)$,
	and pick $\ag_{n+1}\in N$ big enough.
	Then we can (in $V$) fuse this sequence into a $p_\omega\in P$.
	Note that $\dom(p_\omega)\subseteq N\cap I$.
	If $G$ is $P$-generic over $V$ and $p_\omega\in G$, 
	then $p_n\in G$ and $\{p^{[\bar a]}_n:\, \bar a\in\pP(\ag_n)\}$
	is predense below $p_n$, so some $p^{[\bar a]}_n\in G$,
	and $p^{[\bar a]}_n\in D_n\cap N$. 

(4)
is a special case of (5): Set $I_1\DEFEQ I_{<i}$ and
	$I_2\DEFEQ \{i\}$. So $\n\eta_i$ is $V_{<i}$-generic in
$V_{\leq i}$ and therefore in $V_{<\infty}$ as well.

(5)
Set $P\DEFEQ \nwflim_I(Q_i)$, $P_1\DEFEQ \nwflim_{I_1}(Q_i)$, and
$\n P_2$ (the $P_1$-name of) $\nwflim_{I_2}(Q_i)$.

There is a natural map $\phi: p\mapsto (p_1,\n p_2)$ 
from $P$ to $P_1\ast \n P_2$:
$p_1\DEFEQ  p\on I_1$, and $\n p_2$ is defined by
	$\dom(\n p_2)\DEFEQ\dom(p)\setminus I_1$ and
	$\pB {\n p_2} i(\bar\eta)\DEFEQ
	\pB p i ((\n\eta_i)_{(i\in I_1)}^\frown\bar\eta)$.

It is clear that $\phi$ preserves $\leq$.
We claim that it is dense and preserves $\incomp$.
Assume $\phi(p)=(p_1,\n p_2)$, $\phi(q)=(q_1,\n q_2)$, and
$(r_1,\n r_2)\leq (p_1,\n p_2),(q_1,\n q_2)$.
We have to find a $r'\leq_P p,q$ such that $\phi(r')\leq (r_1,\n r_2)$.

$r_1$ forces that $\n p_2$, $\n q_2$ and $\n r_2$
correspond to approximating sequences 
$(\n \ag^p_n)$, $(\n \ag^q_n)$ and $(\n \ag^r_n)$.
As in (1) we can find an $r'_1\leq r_1$ with an approximating
sequence $(\ah_n)$ such that 
$\ah_n$ decides $\n \ag^i_n$ (for $i\in\{p,q,r\}$) 
in a way such that $\n \ag^r_n$ is stronger than both
$\n \ag^p_n$ and $\n \ag^q_n$. Then we can concatenate
$(\ah_n)$ with the $(\n \ag^r_n)$ to an approximating sequence
to some $r'\in P$. Then $r'\leq p,q$ and $\phi(r')\leq (r_1,\n r_2)$.

(6)
By induction on $\epsilon$. The successor step follows from (5).
	Let $\cf(\epsilon)>\omega$. Then 
	the nwf-limit as well as the cs-limit are just the unions of
	the smaller limits, and therefore equal by induction.
	If $\cf(\epsilon)=\omega$, then the nwf-limit as well as the cs-limit 
	are the full inverse limits of the iteration system, and 
	therefore again equal by induction.

(7) follows from (6).
\end{proof}

We will also use the following fact:
\begin{Lem}\label{lem:lem35} Assume that $\n S$ is a $P_{<i}$-name for an 
element of $Q_i$, that
$q\on I_{<i}$ reads $\n S$ continuously and that $q\forc \n\eta_i\in\n
S$. Then $q\on I_{<i}$ forces that $\pB q i (\bar{\n \eta}) \leq_{Q_i} \n S$.
\end{Lem}
\begin{proof}
Assume otherwise. Then there is an approximation $\ag$ of $p\DEFEQ
q\on I_{<i}$, an  $\bar a\in\pP(\ag)$ and a
$t\in\Tmax^i$ such that
$p^{[\bar a]}$ forces $t$ to be in $\pB q i (\bar{\n \eta})$ but not in $\n S$.
Let $\bar a^+$ be $\bar a^\frown t$. Then $\bar a^+$ is a possible 
value of some approximation of $q$, and $q^{[\bar a^+]}$ forces
that $\n\eta_i\notin \n S$, a contradiction.
\end{proof}

\begin{uRem}
  The iteration technique defined here also works for larger classes of
  forcings, e.g., for the tree forcings $\bQ^\text{tree}_0$ of~\cite{MR1613600}
  mentioned already. If we assume additional properties such as bigness and
  halving, we could also use lim-inf forcings. It is also possible to
  extend the construction to non-total orders, or to allow $\Tmax^i,\mu^i$
  to be $P_{<i}$-names. 
\end{uRem}

\section{The ideal $\myI$}\label{sec:myI}

To every tree forcing such as $Q$ defined in Section~\ref{sec:Q} (and many
other tree forcings as well) there is an associated ideal $\strongmyI$ and a
notion of measurability.
We will also use $\myI$, the $<2^\al0$-closure of  $\strongmyI$, and the
associated notion of weak measurability.
The application in this paper of a nw-iteration will be: for certain $Q$ we can
force weak measurability for all definable sets.

\begin{Def}\label{def:IstrongI}
\begin{itemize}
\item
The ideal $\strongmyI$ on the reals is defined by:
$X\in \strongmyI$ if for all $S\in Q$ there is a $T\leq S$
such that $X\cap\lim(T)=\emptyset$.
\item
$\myI$ is the $<2^\al0$-closure of $\strongmyI$.
\item
$X$ has weak measure 1 if $\bR \setminus X \in \myI$.
$X$ has strong measure 1, if $\bR \setminus X \in \strongmyI$.
\end{itemize}
\end{Def}

\begin{uNotes}
\begin{itemize}
\item
Of course these notions depend on the forcing $Q$, so it might be more exact to
use notation such as $\strongmyI_Q$ or $\strongmyI_{(\Tmax,\mu)}$ etc. In this
paper this is not necessary, since we will always use a fixed $Q$.
\item
We use the phrase ``measure 1'' although the ideals $\strongmyI$ and $\myI$ are
not related to a measure (they are not even ccc).
\item
Of course, if CH holds, then $\myI=\strongmyI$.
\item
$\strongmyI$ is always nontrivial (i.e.,  $\lim(\Tmax) \notin \strongmyI$),
but this is not clear for $\myI$.
\end{itemize}
\end{uNotes}

$F: Q\fnto Q$ is a witness for $X\in \strongmyI$
if $F(S)\leq S$ and $X\cap\lim(F(S))=\emptyset$
for all $S\in Q$.

So every $X\in \strongmyI$ is contained in a set $\bigcap
\{\omega\ho\setminus\lim(F(S)):\, S\in Q\}$.\footnote{Note that this is not a
countable intersection.}

\begin{Lem}\label{lem:strongIissigma}
$\strongmyI$ is a non-trivial $\sigma$-ideal.
\end{Lem}
\begin{proof}
This follows from fusion:
Assume $X_i\in \strongmyI$ ($i\in\omega$) and $S=S_0\in Q$.
Pick any front $F_0\in S_0$, so $S_0=\bigcup_{t\in F_0} S_0^{[t]}$.
For each $t\in F_0$ pick an
$S_{1,t}\leq S_0^{[t]}$ such that $\lim(S_{1,t})\cap X_1=\emptyset$.
Set $S_1\DEFEQ \bigcup_{t\in F_0} S_{1,t}$. So $S_1\in Q$,
and $F_0$ is a front in $S_1$.
Pick a $1$-dense  front $F_1$ in $S_1$ (purely) stronger than $F_0$.
Iterate the construction. Fusion produces a
$T< S$ such that $\lim(T)\cap X_i=\emptyset$ for all $i\in \bN$. 
\end{proof}

For example, if  $Q$ is Sacks forcing,
then $\strongmyI$ is called Marczewski ideal.
$X\in\strongmyI$ iff in every 
perfect set $A$ there is a perfect subset $A'$ of $A$ such that
$A'\cap X=\emptyset$. So if $X$ is Borel 
(or if $X$ has the perfect set property, 
e.g., $X$ is $\mS11$),
then $X\in\strongmyI$ iff $X$ is countable.
$\strongmyI$ is not a ccc ideal:
For $A\subseteq \omega$, set 
\[X_A\DEFEQ\{f\in 2\ho:\, (\forall n\notin
A)\,f(n)=0\}.\]
Clearly $X_A\cap X_B=X_{A\cap B}$, and $\card{X_A}=2^{\card{A}}$. So if
$\{A_i:\, i\in 2^\al0\}$ is an almost disjoint family, then 
$\{X_{A_i}\}$ is a family of closed sets not in $\strongmyI$ such that
$X_{A_i}\cap X_{A_j}$ is finite 
for $i\neq j$. 

For a Borel ccc ideal $I$, ``$X\subseteq \bR$ is measurable''
can be defined by ``there is a Borel
set $A$ such that $A\Delta X\in I$''.
(Usually the basis of the ideal is simpler, e.g.,
one can use open sets instead of Borel sets for meager, or $G_\delta$ sets for
Lebesgue-null.) Equivalently, $X$ is measurable iff for every $I$-positive
Borel set $A$ there is an $I$-positive Borel set $B\subseteq A$ such that
either $B\cap X\in I$ or $B\setminus X\in I$.
For non-ccc ideals that do not live on the Borel sets, 
this second notion is usually the one used to define measurability:
\begin{Def}\label{def:measurable}
\begin{itemize}
\item
$X\subseteq \bR$ is measurable if for every 
$T\in Q$ there is an $S\leq_Q T$ such that either 
$\lim(S)\cap X\in \strongmyI$ or $\lim(S)\setminus X\in \strongmyI$.
\item $X\subseteq \bR$ is weakly measurable if for every 
$T\in Q$ there is an $S\leq_Q T$ such that either 
$\lim(S)\cap X\in \myI$ or $\lim(S)\setminus X\in \myI$.
\end{itemize}
\end{Def}

Since $\myI$ is the bigger ideal, measurability implies weak measurability.

In the rest of the paper, we will construct a specific $Q$ and a nwf-iteration
$P$ and show that $P$ forces all definable sets to be weakly measurable:
\begin{Thm}\label{thm:defaremeas}
  Assume CH and that $Q$ satisfies the Ramsey property \ref{def:ramseyprop}.
  Then there is a proper, $\al2$-cc, $\omega\ho$-bounding p.o. $P$ forcing that
  every set of reals which is (first-order) definable using a parameter in
  $L(\bR)$ is weakly measurable.
\end{Thm}
We will see in Lemma~\ref{lem:Qexists} that there is such a $Q$, and the
Theorem will be proven by~\ref{lem:onlyetaalphamatter},~\ref{lem:oneetalalleta}
and~\ref{cor:final}.

\begin{Rem}\label{rem:refq}
  It is natural to ask whether in our forcing extension
  every definable set is measurable (and not
  just weakly measurable, as stated in the theorem).
  This seems unlikely, but it is not clear how to prove it.
  It is not even clear how to prove that in our forcing model
  $\strongmyI\neq \myI$ (i.e., that $\add(\strongmyI)<2^{\al0}$).
  (Of course, $\strongmyI=\myI$ would trivially imply that measurable sets
  and weakly measurable sets are the same, so in particular that all
  definable sets are measurable.)
\end{Rem}

Let us first list some facts about (weak) measurability:

\begin{Lem} Every Borel set is measurable.
The family of measurable sets 
is closed under complements and
countable unions; the same holds for weakly measurable sets.
\end{Lem}

\begin{proof}
Closure under complement is trivial.

Every closed set is measurable: Let $X=\lim(T')$ be closed
and $T\in Q$. If there is a $t\in T\setminus T'$ then 
$S\DEFEQ T^{[t]}$ satisfies $\lim(S)\cap X=\emptyset$.
Otherwise $T\subseteq T'$ and $S\DEFEQ T$ satisfies
$\lim(S)\setminus X=\emptyset$.

Assume that $(X_i)_{i\in\omega}$ is a sequence of weakly measurable sets
and that $T\in Q$. If for some $i\in\omega$ there is an $S\leq T$
such that $\lim(S)\setminus X_i\in\myI$ then the same
obviously holds for $\bigcup_{i\in\omega} X_i$. So assume that
for all $i\in\omega$ and $T'\leq T$ there is an $S\leq T'$ such
that $\lim(S)\cap X_i\in\myI$. 
Now repeat the proof of \ref{lem:strongIissigma}.

The same proof also shows that the measurable sets
are closed under countable unions.
\end{proof}

$\myI$ could be trivial (i.e., $\cov(\strongmyI)$ could be less than
$2^{\al0}$). If $\myI$ is ``everywhere nontrivial'', then $\myI$ and
$\strongmyI$ are the same on measurable (in particular, Borel) sets:

\begin{Lem}\label{lem:IandstrongI}
Assume that $\lim(S)\notin\myI$ for all $S\in Q$.  Then $\myI$ and $\strongmyI$
agree on measurable sets. I.e., if $X$ is measurable and
$X\in \myI$, then $X\in\strongmyI$.
\end{Lem}

\begin{proof}
 For every $T\in Q$ there is an $S\leq_Q T$ such
that $\lim(S)\cap X\in \strongmyI$:  
Otherwise $\lim(S)\setminus
X\in\strongmyI\subseteq \myI$, 
a contradiction to $X\in\myI$ and $\lim(S)\notin\myI$.
So by the definition of $\strongmyI$ there is a $S'\leq_Q S\leq_Q T$ 
such that $ \lim(S')\cap X=\emptyset$. So $X\in\strongmyI$.
\end{proof}

Since any Borel set $B$ is measurable, $B\in\strongmyI$ iff 
$(\forall S\in Q)\lim(S)\not\subseteq B$, so we get:
\begin{Fact}
For a Borel code $B$, the statement ``$B\in\strongmyI$'' is $\mP12$ and therefore
invariant under forcing.
\end{Fact}

On the other hand, since $\strongmyI$ is not a Borel ideal (i.e., not
every $X\in\strongmyI$ is contained in a Borel set $B\in\strongmyI$),
there is no reason why $X\in\strongmyI$ should be upwards absolute
between universes.

For later reference, we will reformulate 
the definition of $\strongmyI$:
If $S\in Q$, $X\subseteq Q$, $T\in X$ and $T'\leq_Q S,T$, then $\lim(T')\cap
(2^\omega\setminus \bigcup_{R\in X}\lim(R))
\subseteq \lim(T')\setminus \lim(T)=\emptyset$.
So we get:
\begin{Lem}\label{lem:predenseiscoi}
If $X\subseteq Q$ is predense then 
$\bigcup_{T\in X}\lim(T)$ is of strong measure 1.
\end{Lem}

\section{An order with many automorphisms}\label{sec:order}

In this section we assume CH. We will construct an order $I$ and define $P$ to
be the nwf-limit of $Q$ along $I$.
$I$ is $\om2$-like,\footnote{$I$ is $\om2$-like if $\card{I_{<i}}<\al2$ for all
$i\in I$ and $\card{I}=\al2$.}
has a cofinal sequence $j_\alpha$ ($\alpha\in\om2$) and
many automorphisms.
We show that these properties imply that $P$ forces the following:
\begin{itemize}
\item $2^\al0=\al2$,
\item $\myI$ is nontrivial (and moreover $\lim(S)\notin \myI$ for all $S\in Q$),
\item for every definable set $X$,
``locally'' either all or no $\eta_{j_\delta}$ are in $X$ and
\item $\{\eta_{j_\delta}:\, \delta\in\om2\}$ 
	is of weak measure 1 in $\{\eta_i:\, i\in I\}$.
\end{itemize}
In the next section it will be shown that 
the set $\{\eta_i:\, i\in I\}$ is of weak measure 1, which will 
finish the proof Theorem~\ref{thm:defaremeas}

First note that for any $I$ with uncountable cofinality, $P$ makes
the old reals null:
\begin{Lem}\label{lem:nooldreal}
  If $I$ has cofinality $\geq \al1$ 
  and $i\in I$ then $\forc_P V_{<i}\cap\lim(\Tmax)\in \strongmyI$.
\end{Lem}
\begin{proof}
	Let $G_P$ be $P$-generic over $V$.
	If $T\in V[G_P]$ then $T\in V_{<j}$ for some $i<j<\infty$
	because of \ref{thm:fusioncor}(\ref{item:nonewreals}).
	So in $V_{\leq j}$ there is an $S<T$ such that
	$\lim(S)\cap V_{<i}=\emptyset$ (in $V_{\leq j}$ and $V[G_P]$ as well,
	according to \ref{lem:treewithoutoldbranch}).
\end{proof}

\begin{Lem}\label{lem:card}
Assume that CH holds and that $I$ is $\om2$-like.
 Then
\begin{enumerate}
\item $P$ has the $\al2$-cc (and therefore preserves all cofinalities).
\item 
 $P_{<i}\forc$ CH
for each $i\in I$ and 
	$P\forc 2^\al0=\al2$.
\end{enumerate}
\end{Lem}

\begin{proof}
(1)
If $\card{I_{<i}}\leq 2^{\al0}$ then $\card{P_{<i}}\leq 2^{\al0}$:
There are at most $\card{I_{<i}}^{\al0}\leq 2^{\al0}$ may countable 
subsets of $\card{I_{<i}}$. For each $p\in P_{<i}$ with a fixed
domain and each $j\in \dom(p)$
there are $2^{\al0}$ many  possibilities for 
$\pD p j$
and $2^{\al0}$ many possibilities for the Borel definition $\pB p j$.

If CH holds, then the usual delta system lemma applies: If $A\subseteq P$ is a
maximal antichain of size $\al2$ then without loss of generality the domains of
$p\in A$ form a delta system
(i.e.,  there is a countable $x\subseteq I$ 
such that
 $\dom(p_1)\cap\dom(p_2)=x$
for all $p_1\neq p_2\in A$).
Since $I$ is $\om2$-like, $x$ cannot be cofinal.
Let $i$ be an upper bound of $x$. Without loss of generality
$p_1\on I_{<i}=p_2\on I_{<i}$
for $p_1\neq p_2\in A$
(since there are only $\al1$ many elements of $P_{<i}$).
But then $p_1\comp p_2$ by Fact~\ref{facts:restriction}.

Proper and $\al2$-cc imply preservation of all 
cofinalities and cardinalities.

(2)
Let $G$ be $P$-generic over $V$.
Then the reals in $V[G]$ are the union of
the reals in $V_{<i}$. Every real in $V_{<i}$ is read
continuously from a condition $p\in G_{<i}$.  There are only
$\card{P_{<i}}=(2^\al0)^V=\al1$ many conditions, and given a condition
there are only $(2^\al0)^V=\al1$ many possibilities to continuously read a real
from the condition. So there are at most $\al1$ many reals in $V_{<i}$.
And $\eta_{i}\notin V_{<i}$, so  in particular $\eta_{i_1}\neq
\eta_{i_2}$ for $i_1\neq i_2$.
\end{proof}

The following is well known:
\begin{Lem}
If CH holds, then there is
an $\al1$ saturated\footnote{A linear order $\tilde I$
is $\al1$ saturated if ``there are no countable gaps'', more exactly: 
\begin{itemize}
  \item $I$ has neither a smallest or a largest element, i.e., no
    $(-\infty,1)$ and no $(1,\infty)$ gaps.
  \item $I$ does not have a cofinal sequence of order type $\omega$ nor
    a coinitial one of order type $\omega^*$, i.e., 
    no $(\omega,\infty)$ and no $(-\infty,\omega^*)$ gaps.
  \item If $A\subset I$ has order type $\omega$ and $c>a$ for all $a\in A$
    ($c>A$ in short)
    then there is a $b<c$ such that $b>A$. I.e., there are no
    $(\omega,1)$ gaps.
  \item Analogously for $B$ of order type $\omega^*$ and $c<B$.
    I.e., no $(1,\omega^*)$ gaps.
  \item If $A$ has order type $\omega$ and  $B$ has order type $\omega^*$
    and $A<B$, then there is an $x\in I$ such that $A<x<B$. I.e., there
    are no $(\omega,\omega^*)$ gaps.
\end{itemize}}
linear order $\tilde I$ of size $\al1$, and
all such orders are isomorphic.
\end{Lem}

\begin{proof}
Induction of length $\om1$: Assume at stage $\alpha$ we 
have a linear order $L_\alpha$ of size $\om1=2^\al0$.
List all the ($\om1$ many) countable gaps
and add points to fill these gaps. At limits, take the union.
Then at stage $\om1$ we get a saturated order.

Uniqueness is proven by the standard back and forth argument.
\end{proof}

\newcommand{\STW}{\mathfrak S}
\begin{Def}
Let $\STW$ be the set of $0<\alpha<\om2$ 
such that $\cf(\alpha)\in\{1,\om1\}$.
Note that $\STW\subseteq\om2$ is stationary.
\end{Def}

We will now define  the order  $I$ along which we iterate.
(We do this assuming CH.)

Given $\tilde I$ as above, let $I$ 
be the following order:
\[\underbrace{\tilde I}_0+
\underbrace{\{j_1\}+\tilde I}_1+\cdots+
\underbrace{\tilde I}_\omega+
\underbrace{\{j_{\omega+1}\}+\tilde I}_{\omega+1}+\cdots+
\underbrace{\{j_{\omega_1}\}+\tilde I}_{\omega_1}+\cdots
\]
So at stages $\alpha\in \STW$, we add an order of the
type $\{c\}+\tilde I$, in other stages we add just $\tilde I$.

\begin{Facts}\label{facts:cofinalsequ}
\begin{itemize}
\item
 $I$ is $\om2$-like,
\item
 $(j_\alpha)_{\alpha \in \STW}$ is an
 increasing (and therefore cofinal) continuous sequence in $I$, and
\item
 every $j_\alpha$ has cofinality $\al1$ in $I$. 
\end{itemize}
\end{Facts}

Continuous 
means that
$j_\delta=\sup(j_\alpha: \, \alpha\in\STW,\alpha < \delta)$
whenever
$\delta=\sup(\STW\cap \delta )\in \STW$ (which is 
equivalent to $\cf(\delta)=\om1$).

\begin{uNote}
We could just as well define $j_\alpha$ for 
$\alpha$ with cofinality $\om1$ only, or for all
$\alpha\in\om2$ (and require continuity for points of cofinality
$\om1$ only). All these versions are equivalent by simple
relabeling, cf.\  the beginning of the proof of \ref{lem:onlyetaalphamatter}.
\end{uNote}

\begin{Def}
We set $Q_i=Q$ for all $i\in I$ and let $P$ be the 
nwf-iteration of $Q_i$ along $I$.

We will use the notation $I_\alpha$, $P_\alpha$, $V_\alpha$ and $\eta_\alpha$
for $I_{<j_\alpha}$, $P_{<j_\alpha}$, $V_{<j_\alpha}$ and $\eta_{j_\alpha}$.
We set $G_\om2$ to be (the name for) the $P$-generic (in previous notation,
$G_{<\infty}$) and $V_\om2$ the generic extension $V[G_\om2]$ (in previous
notation, $V_{<\infty}$).
\end{Def}

\begin{Lem} (CH) Let $S_0\subseteq \STW$ be stationary.
$P$ forces the following:
\begin{enumerate}
\item $\{\eta_{\delta}:\, \delta\in S\}\notin\myI$ for every
	stationary $S\subseteq \STW $, and
\item
	$\{\eta_{\delta}:\, \delta\in S_0\}\cap \lim(T_0)\notin\myI$
	for every $T_0\in Q$.
\end{enumerate}
\end{Lem}

This lemma implies that in $V_\om2$ the assumption of
 Lemma~\ref{lem:IandstrongI} is satisfied (i.e.,  that $\myI$ is 
``everywhere nontrivial''). 
This lemma holds for all $I$ satisfying \ref{facts:cofinalsequ}.

\begin{proof}
(1)
Assume otherwise, i.e.,  there are $P$-names
$\n F_\zeta$ $(\zeta\in\om1)$ for functions from
$Q$ to $Q$ and $\n S$ for a stationary set 
such that $p_0\in P$ forces
\[ \n F_\zeta(T)\leq T
	\text{ and }(\forall \delta\in \n S)\, (\exists \zeta\in\om1)\,(\forall T\in Q)\,
	\eta_{\delta}\notin\lim(\n F_\zeta(T)).
\]

$P$ forces that
for each $\alpha\in \STW $ there is a $\beta \in \STW $ such that
$\n F_\zeta(T)\in Q^{V_{\beta}}$
for all $T\in Q^{V_{\alpha}}$ and  $\zeta\in\om1$.
We need something slightly stronger:
For every name $\n T$ for an element of 
$Q^{V_{\alpha}}$ and $\zeta\in \om1$ there is a maximal antichain 
$A\subset P$ such that for every $q\in A$ there is a
$P$-name $\n{T}'_q$ 
such that $q$ forces $\n F_\zeta(\n T)=\n{T}'_q$ and $q$ continuously 
reads $\n{T}'_q$. So if 
$q\in G_\om2$ and 
$\beta$ is bigger than $\dom(q)$,\footnote{More formally: if
$j_\beta>i$ for all $i\in\dom(q)$.}
then
$V_{\beta}$ not only contains $T'_q=F_\zeta(T)$, but also knows 
that $T'_q$ will be $F_\zeta(T)$ in $V_\om2$.

Define $f^-(\alpha)$ to be the smallest $\beta$
which is bigger than $\dom(q)$ for every $q\in A$,
where $A$ is an antichain for some $\n T$ and $\zeta\in \om1$
as above. $P$ is $\al2$-cc, every $q\in A$ has countable domain,
and there are only $\al1$ many reals in $V_\alpha$.
So $f^-(\alpha)<\om2$,
and we can define $f(\alpha)$
to be the smallest $\beta\in \STW$ that is larger or equal
to $\max(\alpha,f^-(\alpha))$.

If $\cf(\alpha)=\om1$, then $f(\alpha)$ is the supremum of
$\{f(\gamma):\, \gamma\in\STW\cap \alpha\}$, 
since the reals in $V_\alpha$ are the union of
the reals in $V_\gamma$. So $f$ is continuous.

Then $P$ forces the following:
 Since $\n S$ is stationary,
there is a $\beta\in \n S$ such that $f(\beta)=\beta$.
$V_{\beta}$ can calculate every $F_\zeta$, and
$\n F_\zeta''Q$ is dense in $Q$.
Since $\n \eta_{\beta}$ is a $Q$-generic real over 
$V_{\beta}$, there is 
(for every $\zeta\in\om1$) 
a $T\in Q^{V_{\beta}}$ such that $\n \eta_{\beta}\in\lim(\n F_\zeta(T))$,
a contradiction.

(2): We can assume that $T_0\in V$. 
Again, choose names $\n F_\zeta$
as above, and assume that $p_0\in P$ forces that\\
\centerline{$\n F_\zeta(T)\leq T$
	and $(\forall \delta\in S_0)\, (\exists \zeta\in\om1)\,(\forall T\in Q)\,
	\eta_{\delta}\notin\lim(\n F_\zeta(T))\cap\lim(T_0)$.}
Define $f$ as above, so there is a $\beta>\dom(p)$ such that 
$\beta\in S_0$ and $f(\beta)=\beta$. So the same argument 
proves that $p_0$ forces that $\n \eta_{\beta}\notin \lim(T_0)$,
a contradiction.
\end{proof}

We also get the following:

\begin{Lem}\label{lem:onlyetaalphamatter} (CH) 
For every $C\subseteq \om2$ club,
$P$ forces the following:
\[
  \{\n \eta_i:\, i\in I\}\setminus \{\n \eta_\alpha:\, 
  \alpha\in\STW\cap C\}\in\myI.
\]
\end{Lem}
Again, this lemma applies to all $I$ satisfying \ref{facts:cofinalsequ}.

\begin{proof}

We can assume that $C=\om2$, since
we can just relabel the sequence $\{j_\alpha:\, \alpha\in \STW\cap C\}$:
Set $j'_\alpha\DEFEQ j_\beta$, where $\beta$ is the $\alpha$-th element of
$C\cap \STW$. Then $(j'_\alpha)_{\alpha\in\STW}$
satisfies \ref{facts:cofinalsequ} as well.

Recall Definition~\ref{def:spl} of $Q^f_{A^{\psi}_r}$ and $D^\text{spl}_f$ (for $f:\omega \fnto \omega$
increasing and $r\in 2\ho$). Enumerate all 
increasing $f:\omega \fnto \omega$ in $V$ 
as $f_\zeta$ ($\zeta\in\om1$). (CH holds in $V$.)

{\bf Claim:}
In $V$, we can find $P_\alpha$-names $\n T_\alpha^\zeta$ 
($\zeta\in\om1$, $\alpha<\om2$ successor) for elements of $Q$ 
such that the following is forced by $P_\om2$:
\begin{enumerate}
\item
        The set 
	$\{\n T_{\alpha}^\zeta:\, \alpha<\om2\text{ successor}\}\subseteq Q$ is
	dense
	for all $\zeta\in\om1$.
\item
	$\n T_\alpha^\zeta\in D^\text{spl}_{f_\zeta}$ (in $V_\alpha$ or
	equivalently in $V_\om2)$.\footnote{Recall \ref{def:spl} 
			and \ref{lem:dspl}.}
\item
	If $\beta<\alpha$ is a successor, then
	$\n T_\alpha^\zeta$ has no branch in $V_{\beta}$, and
	for all $i<j_\alpha$
	there is a $\zeta_0$ such that 
	$\n T_\alpha^\zeta$ has no branch in $V_{<i}$
	for all $\zeta \geq \zeta_0$.
\end{enumerate}

Proof of the claim:
Pick for all $\alpha+1$ a function $\phi_{\alpha+1}: 
\om1\fnto I_{<j_{\alpha+1}}\setminus  I_{<j_\alpha}$
which is increasing and cofinal. Also pick an
enumeration $(\n S_{\alpha+1})_{\alpha\in \om2}$
such that $\n S_{\alpha}$ is an $P_\alpha$-name 
and $P$ forces that $Q=\{\n S_{\alpha+1}:\, \alpha\in \om2\}$.
(This is possible since $P$ forces that $Q^{V_\om2}=\bigcup Q^{V_\alpha}$,
cf.~\ref{thm:fusioncor}(\ref{item:nonewreals}).)

To find $\n T_\alpha^\zeta$ ($\alpha$ successor)
note that $P_\alpha$ forces that we can perform
the following construction in $V_\alpha$:
First pick an $S'\leq \n S_\alpha$
such that $S'\in D^\text{spl}_{f_\zeta}$  (cf.~\ref{lem:dspl}(\ref{item:dspl6})). $\cf(j_\alpha)=\al1$, so $S'\in V_{<i}$
for some $i<j_\alpha$. Pick some $i'$ bigger than 
$\max(i,\phi_{\alpha}(\zeta))$ and smaller than $j_\alpha$.
There is a real $r\in V_\alpha\setminus V_{<i'}$
(e.g., $\eta_{i'}$). Therefore there is a $T^\zeta_\alpha<S'$
such that $\lim(T^\zeta_\alpha)\cap V_{<i'}=\emptyset$
(in $V_\alpha$ and $V_\om2$ as well, cf.~\ref{lem:treewithoutoldbranch}).
Let $\n T_\alpha^\zeta$ be a $P_\alpha$-name for $T^\zeta_\alpha$.

The $\n T_\alpha^\zeta$ constructed this way satisfy the claim:
(1): $\n T_\alpha^\zeta\leq \n S_\alpha$, (2):
$D^\text{spl}_{f_\zeta}$ is open dense and absolute, (3):
pick $\zeta_0$ such that $\phi_\alpha(\zeta_0)>i$.
This ends the proof of the claim.

From now on assume $G$ is $P$-generic over $V$.
We work in $V_\om2$ and set $T_\alpha^\zeta\DEFEQ \n T_\alpha^\zeta[G]$.
So if $i\in I$ then the sequence
$(T_{\alpha+1}^\zeta)_{j_{\alpha+1}<i,\zeta\in\om1}$ is in $V_{<i}$.

For all $\zeta\in\om1$, $X_\zeta\DEFEQ 
\bigcup_{\alpha+1<\om2} \lim(T_{\alpha+1}^\zeta)$ is of strong measure 1
 (cf.~\ref{lem:predenseiscoi}).
So the set $Y\DEFEQ \bigcap_{\zeta\in\om1} X_\zeta$ is 
of weak measure 1.
It is enough to show that 
\[
  (\{\eta_i:\, i\in I\}\setminus \{\eta_\alpha:\, \alpha\in\STW\})
  \cap Y=\emptyset.
\]
Assume towards a contradiction that some $\eta_i$ is in $Y$ and
$\eta_i\neq \eta_\alpha$ for all $\alpha\in\STW$.

Let $\alpha\in \STW$ be minimal such that $\eta_i\in V_{\alpha}$ (i.e., 
$i<j_{\alpha}$). %, cf.\ Figure~\ref{fig:ordtree}). 
So $\alpha$ is a successor (but not necessarily 
a successor of a $\beta\in\STW$), and $i>j_\beta$ for all
$\beta\in \STW\cap\alpha$.
So according to (3) there is a $\zeta_0$ such that $\eta_i\notin
\lim(T_{\gamma+1}^\zeta)$ for all  $\zeta>\zeta_0$ and all $\gamma+1\geq
\alpha$.

So we know the following:
$\eta_i\in Y$, i.e.,  
\[\eta_i\in\bigcup_{\gamma+1<\om2} \lim(T_{\gamma+1}^\zeta)\quad\text{ for all }
  \zeta\in\om1.\] 
But 
\[\eta_i\notin \bigcup_{\alpha\leq \gamma+1<\om2}
  \lim(T_{\gamma+1}^\zeta)\quad\text{ for all }\zeta\geq \zeta_0.\]
Therefore
\[\eta_i\in\bigcup_{\gamma+1<\alpha} \lim(T_{\gamma+1}^\zeta)\quad\text{ for all }
  \zeta\geq \zeta_0.\] 
Recall that $V_{<i}$ sees the sequence 
$(T_{\gamma+1}^\zeta)_{\gamma+1<\alpha,\zeta\in\om1}$.
So in $V_{<i}$, some $T\in Q$ forces that for all $\zeta>\zeta_0$
there is a successor $\beta(\zeta)<\alpha$ such that
$\n\eta_i\in \lim(T_{\beta(\zeta)}^\zeta)$.
In $V_{<i}$, $T$ has full splitting for some $f_\zeta\in V$,
$\zeta>\zeta_0$ (see \ref{lem:dspl}(\ref{item:dspl4}), \ref{lem:fullsplitinV} and \ref{thm:fusioncor}(\ref{item:omegaomegabounding})).

Let $r$ be a real in 
$V_{<i}\setminus \bigcup_{\gamma+1<\alpha}{V_{\gamma+1}}$.
Pick in $V_{<i}$ a $T'\leq T$ such that $T'\in Q^{f_\zeta}_{A^{\psi}_r}$ 
(cf.~\ref{lem:dspl}(\ref{item:dspl5}))
and $T'$ decides $\beta(\zeta)$.
Then $T'$ forces that $\n\eta_i\in \lim(T'\cap T_{\beta(\zeta)}^\zeta)$,
a contradiction to $T'\incomp T_{\beta(\zeta)}^\zeta\in V_{<i}$
(because of (2), either $T_{\beta(\zeta)}^\zeta$ is in $Q^{f_\zeta}_{A^{\psi}_s}$
for some old real $s$, or incompatible to all $Q^{f_\zeta}_{A^{\psi}_s}$).
\end{proof}

\begin{figure}[tb]
\begin{center}
\scalebox{0.65}{\input{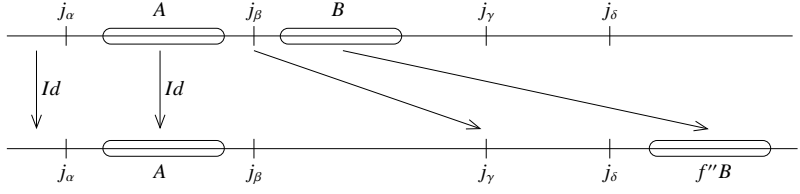}}
\end{center}
\caption{\label{fig:auto}An automorphism $f$.}
\end{figure}

We call $f$ an automorphism if it is a $<$-preserving
bijection from $I$ to $I$.

If $f: I\fnto I$ is an automorphism, then $f$ defines an automorphism of $P$
in a natural way as well (provided of course that $f(i)=j$ implies $Q_i=Q_j$,
but in our case all the $Q_i$
are the same). Also, $f$ defines a map on all $P$-names, and we have:
$p\forc \varphi(\n \tau)$ iff $f(p)\forc \varphi(f\n\tau)$. 

If $\forc_P\n x\in V_{<i}$, then there is a $V_{<i}$-name $\n \tau $ 
such that $\forc_P \n x=\n\tau$.
If $f\on I_{<i}$ is the identity, then $f(\n \tau)=\n \tau$.
So in this case $p\forc \phi(\n \tau)$ iff $f(p)\forc \phi(\n\tau)$.
Also, if $f\on \dom(p)\cap I_{<i}$
is the identity
then $\pB p i(\bar \eta)=\pB{f(p)}{f(i)}(\bar \eta)$.

\begin{Lem}\label{lem:manyaut}
The following holds for $I$ (see Figure~\ref{fig:auto}):
If $\alpha<\beta<\gamma<\delta$ are in $\STW $, and if $A\subseteq I_{\beta}$
and $B\subseteq I\setminus I_{\beta}$  are countable,
then there is an automorphism $f$
of $I$ 
such that $f\on (I_\alpha\cup A)$ is the identity,
$f(j_\beta)=j_\gamma$ and $f'' B>j_{\delta}$.
\end{Lem}

\begin{proof}
For every $i<j\in I$, $I_{<i}$ and $\{k:\, i<k<j\}$ are 
isomorphic and also isomorphic to $\tilde I$ (since they are all $\al1$
saturated
linear orders  of size $\al1$). If
$A\subset I$ is countable, then there are $i<A<j$, and for all
such $i,j$ the sets $\{k:\, i<k<A\}$ and $\{k:\, A<k<j\}$ are again
isomorphic to $\tilde I$.
 Also, $I_{>i}$ is isomorphic to $I$ (since
$\om2\setminus \alpha$ is isomorphic to $\om2$).

So assume $\alpha<\beta<\gamma\in \STW$, $A<i<j_\beta$ countable, $i>j_\alpha$.
Then $I_{<j_\beta}\setminus I_{\leq i}\cong I_{<j_\gamma} \setminus I_{\leq i}\cong \tilde I$.
Also, if $B\subset I$ is countable, $\delta\in \STW$ and
$B>j_\beta$, then there is an $j_\beta<i<B$,
and $I_{<i}\cong I_{<j_\delta}\cong \tilde I$, $I\setminus I_{\leq i}\cong 
I\setminus
I_{\leq j_\delta}\cong I$. Now combine these automorphisms.
\end{proof}

\begin{Lem}\label{lem:oneetalalleta}
For $\beta\in\om2$ set $Y_\beta\DEFEQ \{\n \eta_{\gamma}:\, \gamma\in\STW,\,\gamma>\beta\}$.
$P$ forces the following:
If $X$ is a set of reals defined with a parameter 
$x\in \bigcup_{i\in I} V_{<i}$,
and if $T\in Q$, then
there is an $S \leq T$ and a $\beta\in\om2$  such that
either $\lim(S)\cap X\cap Y_\beta=\emptyset$ or 
$(\lim(S)\setminus X)\cap Y_\beta=\emptyset$.
\end{Lem}
This lemma holds for all $I$ satisfying \ref{facts:cofinalsequ}
and \ref{lem:manyaut}.

Note that every real in $V_\om2$ is in $\bigcup_{i\in I} V_{<i}$.

We will see in the next section that (using additional assumptions)
$Y_\beta$ is a weak measure 1 set. Then this lemma implies that
$X$ is weakly measurable, i.e., Theorem~\ref{thm:defaremeas}. 
Because of \ref{lem:onlyetaalphamatter},
it will be enough to show that 
the set $\{\eta_i:\, i\in I\}$ is of weak measure 1.

\begin{proof}
Assume $\n X=\{r:\, \varphi(r,\n x)\}$ and fix some $\n T$.
Some $p_0$ forces that $\n x$ and $\n T$ are in $V_{\alpha}$, so
without loss of generality
${\n x}, \n T$ are $P_\alpha$ names and $\dom(p_0)\subset I_{\alpha}$.
Pick a $p_1\leq p_0$, $p_1\in  P_\alpha$ such that $p_1$ continuously
reads $\n T$.
Fix some $\beta>\alpha$.
Then $p_2\DEFEQ p_1\cup \{(j_\beta,\n T)\}$
is an element of $P_{\leq j_\beta}$ (since $\n T$ is read continuously).

Let $p\leq p_2$ decide $\varphi(\n\eta_{{\beta}},\n x)$. 
Without loss of generality
$p\forc \varphi(\n\eta_{{\beta}},\n x)$.
$p\on I_\beta$ 
forces that $\n S\DEFEQ \pB p {j_\beta}(\bar{\n \eta})\leq_Q \n T$
(since $p\leq p_2$).

Assume towards a contradiction that for some $q\leq p$, $\gamma\in \STW $ and
$\gamma>\beta$
\[q\forc \n\eta_{\gamma}\in\lim (\n S)\ \& \ \lnot \varphi(\n
\eta_{\gamma},\n x).\]
Note that $q\on I_\gamma$ reads $\n S$ continuously and forces
that $\pB q {j_\gamma} (\bar{\n \eta})\leq_Q \n S$
(cf.~\ref{lem:lem35}). 

Set $A\DEFEQ \dom(p)\cap I_\beta$ and
$B\DEFEQ \dom(p)\cap I_{>j_\beta}$.
Let $j_{\delta}$ be bigger than $\dom(q)$, and 
let $f$ be an automorphism of $I$ such that
$f\on (I_\alpha\cup A)$ is the identity,
$f(j_\beta) = j_{\gamma}$ and $f'' B>\dom(q)$
(cf.~\ref{lem:manyaut} or Figure~\ref{fig:auto}).

$\dom(f(p))\cap \dom(q)\subseteq A\cup \{j_\gamma\}$.
$f(p)\on A=p\on A\geq q\on A$,
and $q\on I_\gamma$ forces that
\[
  \pB {f(p)} {j_{\gamma}}(\bar{\n \eta})=\pB p {j_\beta}(\bar{\n \eta})=\n S
  \geq_Q \pB q {j_\gamma} (\bar{\n \eta}).
\]
So $f(p)$ and $q$ are compatible, a contradiction to
$f(p)\forc \varphi(\n\eta_{{\gamma}},\n x)$.
\end{proof}

\section{A very non-homogeneous tree}\label{sec:tree}

For the proof of Theorem~\ref{thm:defaremeas} it remains to be shown that
$\{\eta_i:\, i\in I \}$ is of weak measure 1.
For this we will need a certain Ramsey property for $Q$.  

\begin{Def}
A subtree $T$ of $\Tmax$ is called $(n,r)$-meager 
if
$\mu_T(t)<r$
for all $t\in T$ with length at least $n$. 
\end{Def}

\begin{Lem}
If $T$ is meager for some $(n,r)$,  then $\lim(T)\in \strongmyI$.
\end{Lem}
\begin{proof}
For any $S\in Q$ there is an $s\in S$ of length at least $n$ such that
$\mu_S(s)>r$. So there is an immediate successor $t$ of $s$ in $S$ such that
$t\notin T$.  Then $\lim(S^{[t]})\cap \lim(T)=\emptyset$.
\end{proof}

\begin{Def}
Let $M,N$ be natural numbers. $N\rightarrow M$ means: If
\begin{itemize}
\item $r_1,\dots,r_M\in \Tmax$ such that $\length(r_i)>N$,
\item $t\in  \Tmax$ such that $r_i\incomp t$ for $1\leq i\leq M$,
\item $A\subseteq \SUCC(t)$ such that $\mu(A)>N$,
\item $f_i:A\fnto \Tmax^{[r_i]}$ for $1\leq i\leq M$,
\end{itemize}
then there is a $B\subseteq A$ such that
\begin{itemize}
\item $\mu(B)>M$ and
\item $\{s\in \Tmax:\, (\exists i\leq M)\, (\exists t\in B)\, s\preceq f_i(t)\}$ is $(N,1/M)$-meager. 
\end{itemize}
\end{Def}

\begin{Def}\label{def:ramseyprop}
A lim-sup tree-forcing $Q$ is strongly non-homogeneous if
$\mu$ is sub-additive\footnote{$\mu(A\cup B)\leq \mu(A)+\mu(B)$.}
and for all $M$ there is an $N$ such that $N\rightarrow M$.
\end{Def}
There are many similar notions of bigness, see e.g., \cite[2.2]{MR1613600}.

\begin{Lem}\label{lem:Qexists}
There is a forcing $Q$ that is strongly non-homogeneous.
\end{Lem}

\begin{proof}
First note that it is enough to show that
for each $M$ there is an $N$ such that $N\rightarrow^- M$,
where $N\rightarrow^- M$ is defined as above but
with just one $r$ and $f$ instead of $M$ many.
To see this, just set $K_0\DEFEQ M^2$ and 
find $K_i$ such that $K_{i+1}\rightarrow^- K_i$.
Then $K_M \rightarrow M$. (Here we use that $\mu$ 
is sub-additive, since we need that the union of $m$ 
many $(n,x)$-meager trees is $(n,x\cdot m)$-meager.)

We will construct $\Tmax$ and $\mu$ by induction.
We define $s\lhd t$ by: $\length(s)<\length(t)$
or $\length(s)=\length(t)$ and $s$ is lexicographically smaller
than $t$. 

Fix some $t\in\omega\hko$.  Assume that we already decided which $s\lhd t$ will
be elements of $\Tmax$ and that we already defined the set of successors of all
these $s$ as well as the measure of their subsets. Assume that we have decided
to put $t$ into $\Tmax$. So we have to define $\SUCC(t)$ and the measure on it.

Let $m_t$ be the number of nodes $s\lhd t$ already defined, including the
already defined successors of $s$ for $s\lhd t$.
Set $M_t\DEFEQ (2 m_t)^{m_t}$ . 
Then we define $\SUCC(t)$ to be of size ${M_t}^{m_t}$.\footnote{
We can e.g., set $\SUCC(t)\DEFEQ \{t^\frown k:\, 0\leq k<M_t^{m_t}\}$.}
For $A\subseteq \SUCC(t)$ we set $\mu(A)\DEFEQ
\log_{M_t}((\card{A}/{M_t})+1)$.

Then $0\leq \mu(A)<m_t$,
$\mu(A)=0$ iff $A=\emptyset$, and $\mu$ is strictly monotonous and
sub-additive.\footnote{Since the function $g(x)\DEFEQ\log_{M_t}(x+1)$ 
is concave and satisfies $g(0)=0$.} If $A,B\subseteq \SUCC(t)$ and 
$\card{B}\geq \card{A}/{M_t}$, then $\mu(B)>\mu(A)-1$.
If $\card{B}\leq m_t$ then $\mu(B)<1/{m_t}$.
If $\mu(\SUCC(t))>M$, then $m_t>M$.

Now fix an arbitrary $M\in\omega$. There is an $N_0$ such that 
$\mu(A)<1/M$
for all $s$ with $\length(s)>N_0$ and all $A\subseteq \SUCC(s)$
with $\card{A}<m_s$.
(Just 
note that $m_s$ strictly increases with $\length(s)$.)
Let $N$ be larger than $M+1$ and $N_0$.

So assume that $r\incomp t\in \Tmax$,
$\length(r)>N\geq N_0$, $A\subseteq \SUCC(t)$,
$\mu(A)>N\geq M+1$ (in particular $m_t>M$),
and $f: A\fnto \Tmax^{[r]}$.

Set $X\DEFEQ\{s'\succeq r:\, s'\lhd t,\, \length(s')\geq N\}$. 
Enumerate $X$ as $\{s_0,\dots,s_{l-1}\}$ (for some $l\geq 0$).
Set $A_0\DEFEQ A$. Assume that $A_n$ is already defined, and
define
\[S_n\DEFEQ \{s'\in \Tmax:\, (\exists t'\in A_n)\, s'\preceq f(t')\}.\] 
If $n>0$ assume that $\card{\SUCC_{S_n}(s_{n-1})}\leq 1$ and that
$\card{A_n}>\card{A_{n-1}}/(2 m_t)$.

Then we define $A_{n+1}$ as follows:
Since $s_n\in X$, $\card{\SUCC(s_n)}<m_t$. By a simple pigeon-hole
argument, there is an 
 $A_{n+1}\subseteq A_n$ such that $\card{A_{n+1}}>\card{A_n}/(2 m_t)$ 
and $\card{\SUCC_{S_{n+1}}(s_{n})}\leq 1$. So in the 
end we get a $B\DEFEQ A_l$ with cardinality at least 
$\card{A}/(2 m_t)^{m_t}=\card{A}/M_t$, i.e.,  
$\mu(B)>\mu(A)-1\geq M$.
Also, $\card{\SUCC_{S_l}(s')}\leq 1$ for every $s'\in X$,
so $\mu_{S_{l}}(s')\leq 1/M$ (since $\length(s')$
was sufficiently large).

We claim that $B$ is as required.
We have to show that $S_l$ is $(N,1/M)$-meager.
Pick an $s'\in S_l$ of length $\geq N$. 
We already dealt with the case $s'\in X$.
Otherwise $s'\rhd t $ (note that $s'\neq t$ since $s'\incomp t$).
In this case
$\card{\SUCC_{S_l}(s')}\leq \card{\SUCC_{\Tmax}(t)}\leq m_{s'}$.
So $\mu(\SUCC_{S_l}(s'))\leq 1/M$, since $\length(s')>N_0$. 
\end{proof}

\begin{Lem}\label{lem:strongnonhom}
If $Q$ is strongly non-homogeneous, then $P$ forces the following:
If $r\in\lim(\Tmax)\setminus  \{\n\eta_i:\, i\in I\}$
then there is a $T\in V$ such that $r\in\lim(T)$ and 
$T$ is $(1,1)$-meager.
\end{Lem}

If additionally the assumptions of Lemma~\ref{lem:card} hold, then 
there are only $\al1$ many $T\in V$, and 
$\al1<(2^\al0)^{V_\om2}$. This implies 
that the set $\{\n\eta_i:\, i\in I\}$ 
is of weak measure 1:\\
If $r\in\lim(\Tmax)\setminus  \{\eta_i:\, i\in I\}$,
then $r\in \bigcup_{T\in V\text{ meager}}\lim(T)\in \myI$.

\begin{proof}
Fix a $P$-name $\n r$ for a real and a  $p\in P$
such that $p\forc \n r\notin \{\n\eta_i:\, i\in I\}$.
We will show that there is a
$p_\omega\leq p$ and a $(1,1)$-meager tree $T$ such 
that $p_\omega\forc \n r\in\lim(T)$.

We will by induction construct $p_n\in P$,
approximations $\ag_n$, $k_n\in \omega$ and $i_n\in u_n=\dom(\ag_n)$ such that
\begin{enumerate}
\item\label{item:jkhwe_1} $p_{n+1}\leq_{\ag_n} p_n$,
$\ag_{n+1}$ is purely stronger than $\ag_n$. 
\item\label{item:jkhwe_2} $\ag_n$ is $n$-dense at $i_n$.
\item\label{item:jkhwe_6} the sequence $(i_n)_{n\in\omega}$
	covers $\bigcup \dom(p_n)$ infinitely often.
\item\label{item:jkhwe_vier} $k_n \rightarrow \max(n+1,\card{\pP(\ag_n)})$.
\item\label{item:meager} If $n>0$, then for
        each $\bar a\in \pP(\ag_n)$, $p_n^{[\bar a]}$ 
	forces a value to $\n r\on k_n$, 
	and the tree $\{r^{\bar a}_n:\, \bar a\in\pP(\ag_n)\}\subseteq
        \Tmax\on k_n$ is $(k_{n-1},1)$-meager.
\end{enumerate}
(1)--(3) allow us to fuse the $(p_n)_{n\in\omega}$ into a $p_\omega\leq p$
(cf.~\ref{lem:fusion}), and (5) implies that the tree of all initial segments
of $r$ compatible with $p_\omega$ is meager. 

We start by picking any $i_0\in\dom(p)$, some $p$-approximation $\ag_0$ that is
$0$-dense at $i_0$, a $k_0$ satisfying (\ref{item:jkhwe_vier}).  So assume by
induction we have found $p_n$, $\ag_n$ and $k_n$ satisfying (1,2,4). 
\begin{itemize}
\item[(a)]
  Set  $p\DEFEQ p_n$, $\ag\DEFEQ \ag_n$, $M\DEFEQ \card{\pP(\ag_n)}$
  and $N\DEFEQ k_n$. So we have $N\rightarrow M$.
\item[(b)]
	Choose the position $i_{n+1}\in\dom(p_n)$ 
	according to some simple bookkeeping.
	This takes care of (\ref{item:jkhwe_6}).
	Set $j\DEFEQ i_{n+1}$. 
\item[(c)]
	Find a $p_1\leq_{\ag} p$ and $m>N$ such 
	that $p_1\forc(\n\eta_j\on m\neq \n r\on m)$
	and for all $\bar a\in \pP(\ag)$ the condition
        $p_1^{[\bar a]}$ determines
	$\n\eta_j\on m$ and $\n r\on m$.

        (How to do this?
	First apply pure decision
	\ref{cor:restrtoa}(\ref{item:puredecision}) to 
	get a $p'\leq_{\ag} p$ such that for all 
	$\bar a\in\pP(\ag)$ there is an $m^{\bar a}>N$
	and $\eta^*\neq r^*$
	such that $p'^{[\bar a]}\forc (r^*=\n r\on m^{\bar a},
	\eta^*=\n\eta_j\on m^{\bar a})$.
	Then we apply pure decision again to get $p_1\leq_{\ag} p'$ 
	determining $\n r$ and $\n \eta_j$ up to 
        $\max\{m^{\bar a}:\, \bar a\in\pP(\ag)\}$.)
\item[(d)]
	Pick a $p_1$-approximation $\ah_1$ which is $\max(n,N)$-dense at $j$
	and (purely) stronger than $\ag$.
\item[(e)]  Pick a $k_{n+1}>m$
       such that $k_{n+1}\rightarrow \max(n+2,\card{\pP(\ah_1)})$.
\item[(f)]	Pick a $q\leq_{\ah_1} p_1$ such that
	$q^{[\bar b]}$ determines 
        %$\n\eta_j\on k_{n+1}$ and 
        $\n r\on k_{n+1}$ up to $k_{n+1}$
	for all $\bar b\in\pP(\ah_1)$.
\end{itemize}
	So far we have taken care of (1--4):
        $q\leq_{\ag}p$,  $\ah_1$ approximates $q$ and 
	witnesses $N$-density (at $j$).
	However, the tree of possible values for $\n r$
	could be very thick in the levels between $k_n$ and $k_{n+1}$.
	We will
	thin out the approximation $\ah_1$ so that
	we still have $(n+1)$-density, and the tree of
	possible values for $\n r$ gets sufficiently thin.
	We do this in two steps:
\begin{figure}[tb]
\begin{center}
\scalebox{0.44}{\input{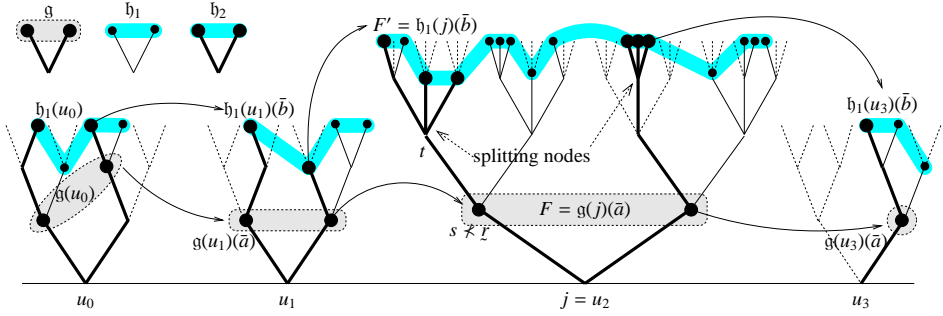}}
\end{center}
\caption{\label{fig:ramsey}$\ah_2$ (bold) is a subapprox.\ of $\ah_1$ and still purely stronger than $\ag$. Here, we assume $\dom(\ah_1)=\{u_0,\dots,u_3\}$,
$j=u_2$, $\bar b\in\pP(\ah_2)$ and $\bar a=\bar b\on\ag$. }
\end{figure}
\begin{itemize}
\item[(g)]
	Find a sub-approximation $\ah_2$ of $\ah_1$
	that is still purely stronger than $\ag$	
	and has only as many splittings as $\ag$,
	apart from one additional split (for each possibility) 
        that witnesses 
	$N$-density at $j$ (see Figure~\ref{fig:ramsey}).

	In more detail: we construct $\ah_2$ the
	following way: Given $\bar b \in \pP_{<i}(\ah_2)$,
        set $\bar a=\bar b\on \ag$.
	We have to define $\ah_2(i)(\bar b)$.
	If $i\neq j$, pick for each $t\in \ag(i)(\bar a)$
	exactly one successor $s\in \ah_1(i)(\bar b)$. So 
	$\ah_2$ makes the branches of $\ag$ longer,
	but does not add any splittings.
	At $j$, we have the front $F\DEFEQ\ag(j)(\bar a)$ 
	and the purely stronger  $n'$-dense front $F'\DEFEQ \ah_1(j)(\bar b)$.
	Recall that $T\DEFEQ \Tcldn^{F'}=\{s:\, s\preceq F'\}$ 
	is the finite tree corresponding to the front $F'$.
	We continue each $t\in F$ in $T$
	uniquely (without splits) until we reach a node $t$ with many
	(i.e.,  $n'$-dense) splittings. We call $t$ ``splitting node''.
	We take all the immediate successors of the splitting node
	and continue them
	uniquely in $T$ until we reach a leaf of $T$, i.e., 
	an element of $F'$. This process leads to a subset
	$F''$ of $F'$. Set
	$\ah_2(j)(\bar b)\DEFEQ F''$.
\item[(h)]
	So we get:
	There are $\card{\pP_{\leq j}(\ag)}\leq  M $ many 
	pairs $(\bar b,t)$, where 
	$\bar b\in \pP_{< j}(\ah_2)$  and $t$ is a splitting node.

	Also, for $\bar b\in \pP_{\leq j}(\ah_2)$,
	there are at most $ M $ continuations 
	of $\bar b$ to some $\bar b'\in \pP(\ah_2)$.

	Such a $\bar b\in \pP_{\leq j}(\ah_2)$
	corresponds to a pair $(\bar a, t)$ as 
	above together with a choice of an (immediate) successor of $t$.
\item[(i)]
	Now we are ready to apply the Ramsey property.
	First fix a $\bar b\in \pP_{<j}(\ah_2)$ and a 
	splitting node $t$. (There are at most $M$ many such pairs.)

	This pair corresponds to a unique $\bar a\in\pP_{\leq j}(\ag)$.
	There are at most $M$ many continuations of $\bar a$ to some
	$\bar c\in\pP(\ag)$. Fix an enumeration $\bar c_1\dots \bar c_M$
	of these possible continuations. According to (c),
        each $\bar c_l$ forces a value
	to $\n r\on m$, call this value $r_l$.

	Back to $\ah_2$. Set $A\DEFEQ \SUCC(t)$ in the tree
	$\Tcldn^{\ah_2(j)(\bar b)}$ (or equivalently
	$\Tcldn^{\ah_1(j)(\bar b)}$).  So $\mu(A)\geq n'>N$.
	For every $s\in A$ there is a unique $s'\succeq s$ such that
	$\bar a \cup \{(j,s')\}\in \pP_{\leq j}(\ah_2)$,
	 and for every
	$s\in A$, $l\in M$ there is a unique 
	$\bar d\in\pP(\ah_2)$ continuing $\bar c_l\in \pP(\ag)$
	and $\bar a\cup\{(j,s')\}$.
	Each such $\bar d$ decides $\n r$ up to $k_{n+1}$. We
        call this value $r^{s,l}$.
	So $r^{s,l}\on m=r_l$.
        According to (d) we know that  $\length(t)>m$, so in
        particular $t\incomp r_l$, according to (c).

	So for every $l\in M$
	we define a function $f_l:A\fnto \Tmax^{[r^l]}\on k_{n+1}$
	by mapping $s$ to  $r^{s,l}$.
	So we can apply the Ramsey property and get a $B\subseteq A$
	such that $\mu(B)>M\geq n+1$, and the tree of possibilities
	for $\n r$ induced by $\bar a,B$ is $(k_n,1/M)$-meager. We repeat that 
	for all pairs $(\bar a,t)$ where 
        $\bar a\in\pP_{<j}(\ah_2)$ and $t$ is a splitting node,
	and get a subapproximation $\ag_{n+1}$ of $\ah_2$ such that
	the tree of possibilities for $\n r$ induced by $\ag_{n+1}$ 
	is $(k_n,1)$-meager (here we again use the
	sub-additivity of $\mu$).
\end{itemize}
     	This results in a sub-approximation $\ag_{n+1}$ of 
	$\ah_2$ (and therefore $\ah_1$) which is still purely
	stronger than $\ag=\ag_n$. Since $\ag_{n+1}$ is 
	a sub-approximation of $\ah_1$, $\card{\pP(\ag_{n+1})}\leq
	\card{\pP(\ah_1)}$, and therefore $k_{n+1}$, $\ag_{n+1}$ satisfy
	(4).
\end{proof}

Note that we did not use the $j_\alpha$ or automorphisms of $I$,
the proof works for all $I$.
In particular, for $I=\{i\}$ we get: 
If $G$ is $Q$-generic over $V$, and if $r\neq \eta$
in $V[G]$, then there is a $(1,1)$-meager $T$ in $V$
such that $r\in \lim(T)$. In particular, such an $r$ cannot 
be $Q$-generic over $V$. So we get:
\begin{Cor}\label{cor:jesh}
  If $Q$ is strongly non-homogeneous then $Q$
  forces that $\n \eta$ is the only $Q$-generic real over $V$ in $V[G_Q]$.
\end{Cor}

\begin{uRem} A similar forcing $Q^\text{JeSh}$ (finitely splitting,  rapidly
increasing number of successors) was used in \cite{MR1398120} to construct a
complete Boolean algebra without proper atomless complete subalgebra.
$Q^\text{JeSh}$ can also be written as lim-sup forcing.  However, the
difference is that the norm in $Q^\text{JeSh}$ is ``binary'' (as e.g., Sacks):
either $s$ has a minimum number of successors, then the norm is large, or the
norm is $0$. Such a norm cannot satisfy a Ramsey property as the one above. For
$Q^\text{JeSh}$ we can only prove Corollary~\ref{cor:jesh} for the ``single
step iteration'', but not Lemma~\ref{lem:strongnonhom} for the iteration.
\end{uRem}
 
We have already mentioned another corollary:
\begin{Cor}\label{cor:final}
If $Q$ is strongly non-homogeneous,
then $P$ forces that
$\{\n \eta_i:\, i\in I\}$ is of weak measure 1.
\end{Cor}

This, together with \ref{lem:onlyetaalphamatter} and \ref{lem:oneetalalleta}
proves Theorem~\ref{thm:defaremeas}.

\begin{uRem}
There are various ways to extend the constructions in this paper.  As already
mentioned, we could use non-total orders $I$ or allow $Q_i$ to be a
$P_{<i}$-name. A more difficult change would be to use lim-inf trees instead of
lim-sup trees. In this case we need additional assumptions such as bigness and
halving.  This could allow us to apply Saccharinity to a ccc ideal
$\strongmyI$, i.e., to force (without inaccessible or amalgamation)
weak measurability of all definable sets.
\end{uRem}

\section{The Cohen model}\label{sec:cohen}
We thank the referee for providing this section.

There is a well known and much simpler way to force that every definable set is
even measurable (not just weakly measurable) with respect to many tree
forcings: Just add many Cohen reals.

Let $\mathbb C^\kappa$ be the forcing notion adding $\kappa$ many Cohen reals
(in a finite support product, or, equivalently, a finite support iteration).
Any $\kappa$ with uncountable cofinality will work.
We call the forcing extension the ``Cohen model''.  If in the ground model
$\kappa^{\al0}=\kappa$, then the continuum has size $\kappa$ in the Cohen
model.

\begin{Lem}
  In the Cohen model, every definable (e.g., projective) set is $Q$-measurable.
\end{Lem}
This works for all $Q$ as in Section~\ref{sec:Q}, in particular for Sacks
forcing, and also many other tree forcings, such as Silver forcing (as was
shown in~\cite{MR2127234}).
So in particular, in the Cohen model all definable sets are Marczewski
measurable (corresponding to $Q$ = Sacks) and have the doughnut property
(corresponding to $Q$ = Silver).

\begin{proof}
  This is similar to, but simpler than, Solovay's argument that all definable
  sets are Lebesgue measurable in the Solovay model.

  Assume that in the Cohen model the parameter $p$ is in the
  union of the intermediate extensions (i.e., already added by
  the first $\alpha$ Cohen reals for some $\alpha<\kappa$) and that
  \[
    X=\{ x:\, \varphi(x,p)\}.
  \]
  for some first order formula $\varphi$. Pick $T\in Q$.
  We can assume without loss of generality (by factoring $\mathbb C^\kappa$)
  that $p$ and $T$ are in $V$.

  Work in $V$ and consider the 
  (countable) forcing notion $T$ (ordered by $\leq_T$, the standard tree
  order). This forcing (which is obviously equivalent to a single Cohen forcing) 
  adds a real $\n c$ that is Cohen over $V$ in the natural topology of $\lim(T)$
  (we call such a real $T$-Cohen, for short).
  In the same way as for ``standard Cohen'' forcing, one can see that
  $\n c$ determines the $T$-generic filter, and $c^*$ is $T$-Cohen
  iff $c^*$ is $\n c[G]$ for some $T$-generic $G$ over $V$.

  In particular, whenever $R$ is some forcing notion, $G_R$ is $R$-generic over $V$ and
  $c^*\in V[G_R]$ is $T$-Cohen (over $V$), then we can factor the extension
  by first adding the $T$-generic $c^*$ and then forcing with some quotient forcing
  to extend $V[c^*]$.
  If $R$ is $\mathbb C^\kappa$, then the quotient forcing is again equivalent
  to $\mathbb C^\kappa$.

  Let $c^*$ be $T$-Cohen (i.e., $T$-generic) over $V$.
  In $V[c^*]$ consider the forcing notion $\mathbb C^\kappa$.
  Since this forcing is homogeneous, 
  either $\forc_{\mathbb C^\kappa}\varphi(c^*,p)$
  or $\forc_{\mathbb C^\kappa}\lnot\varphi(c^*,p)$.
  Without loss of generality assume the former. 
  So in $V$ we can pick some condition $t^*\in T$ such that
  \[
    t^*\forc_T \forc_{\mathbb C^\kappa}\varphi(\n c,p).
  \]

  Let $c^*$ in a $\mathbb C^\kappa$-extension $V'$ of $V$ be any $T$-Cohen real extending $t^*$.
  As described above, we can get $V'$ by first extending $V$ with
  the $T$-generic $c^*$ and then some $\mathbb C^\kappa$-extension of $V[c^*]$.
  In particular, $\varphi(c^*,p)$ holds in $V'$. To summarize:
  \begin{equation}\label{eq:di}
    \text{In the Cohen model $V'$, all $T$-Cohen reals $c^*$ that extend $t^*$
      satisfy $\varphi(c^*,p)$.}
  \end{equation}
  Back in $V$, let $T'$ be the tree $T^{[t^*]}$. So $T'\in Q^V$.
  Set
  \[
    P=\{(t,n):\, n\ge \length(t^*), t\text{ is a subtree of } T'\text{, each maximal branch has height }n\}
  \]
  ordered by end-extension (more exactly: $(t,n)$ is stronger than $(s,m)$ 
  iff $n\ge m$ and $t$ end-extends $s$). 
  Obviously $P$ is equivalent to Cohen forcing as well, and $P$ adds a 
  generic subtree $S$ of $T'$ (and $S$ determines the generic filter). 
  By density, the lim-sup condition
  will be satisfied, so $S$ is in $Q^{V[S]}$.
  In any forcing extension $V'$ of $V[S]$, we get:
  \begin{equation}\label{eq:wijtt}
    \text{Every branch $\nu\in\lim(S)$ is $T$-Cohen over $V$ and extends $t^*$.}
  \end{equation}
  To see this, fix some nowhere dense set $N$ in $V$. Without loss
  of generality $N$ is closed, i.e., corresponds to a nowhere dense subtree $N'$ of $T$. 
  Then (by a simple density argument) 
  there is some $(t,n)$ in the $P$-generic such that each maximal branch
  of $t$ is not in $N'$. So any  $\nu\in \lim(S)$ extends one of the maximal branches
  of $t$, and therefore is not in $N$.

  Now we can finally fix a $\mathbb C^\kappa$-extension $V'$ of $V$.
  We can use the equivalence of $\mathbb C^\kappa$ and $P\ast \mathbb C^\kappa$
  to get in $V'$ some $S\leq_Q T$ such that~\eqref{eq:wijtt} holds.
  Then by~\eqref{eq:di} we get that each $c^*\in\lim(S)$
  satisfies $\varphi(c^*,p)$, i.e., that $\lim(S)\subseteq X$.
\end{proof}

What is the difference between the Cohen model and the model obtained in the
non-wellfounded iteration (let us call it nw-model, for short)? Note that in
our nw-model, the continuum has size $\al2$ (of course we can get larger
continuum as well). One obvious difference is that in the nw-model $\myI$ (the
$<\al2$-closure of $\strongmyI$) is non-trivial (or, in the language of
cardinal characteristics, $\cov(\strongmyI)=\al2$), which is
not the case in the Cohen model for $\kappa\geq \al2$:
\begin{Lem}
  In the Cohen model, $\cov(\strongmyI)=\om1$.
\end{Lem}
\begin{proof}
  The Cohen model is obtained by a finite support product of
  $\kappa$ many Cohen reals. 
  We can write $\kappa$ as the strictly increasing union 
  $\bigcup_{\alpha\in\om1} A_\alpha$ (each $A_\alpha$ of size $\kappa$).
  Let $\mathbb C_\alpha$ be the complete subforcing 
  of $\mathbb C^\kappa$ consisting of the conditions that only use coordinates in
$A_\alpha$.
  Let $G$ be $\mathbb C^\kappa$-generic over $V$, and let 
  $G_\alpha$ be the induced $\mathbb C_\alpha$-generic filters over $V$.
  Then we get:
  \begin{enumerate}
    \item $V[G_\alpha]\cap \omega^\omega$ is a proper
      subset of $V[G_{\alpha+1}]\cap \omega^\omega$.
    \item $V[G]\cap \omega^\omega=\bigcup_{\alpha\in\om1} V[G_\alpha]\cap
      \omega^\omega$.
  \end{enumerate}
  From (1) and Lemma~\ref{lem:treewithoutoldbranch} we know
  that each $V[G_\alpha]\cap \omega^\omega$ is $Q$-null
  in the final Cohen extension;
  so by (2) $\omega^\omega$ is the 
  union of $\al1$ many $Q$-null sets.
\end{proof}
(This argument works not only for the Cohen extension, but also for the random
model and similarly for finite support iteration of Suslin ccc forcings of
length $\al2$; also, it works for other ideals than the ones defined by lim-sup
tree forcings.)

\bibliographystyle{amsplain}
\bibliography{saccharinity}

\end{document}